\documentclass{amsart}
\usepackage[utf8]{inputenc}

\usepackage[english]{babel}

\usepackage{amssymb}
\usepackage{amsmath}
 \usepackage{amsmath,amscd}
\usepackage{amsthm}
\usepackage{enumerate}
\usepackage{cite}
\usepackage{xcolor}
\usepackage[final]{listings}
\usepackage[section]{algorithm}

\makeindex

\title[On the identifiability of ternary forms]
{On the identifiability of  ternary forms}
\date{}

\newcommand{\C}{\mathbb{C}}
\newcommand{\Z}{\mathbb{Z}}
\newcommand{\Pj}{\mathbb{P}}

\newcommand{\N}{\mathbb{N}}

\newcommand{\vect}[1]{\mathbf{#1}}

\newcommand{\imm}{\operatorname{im}}

\newtheorem{thm0}{Theorem}[section]
\newtheorem{prop0}[thm0]{Proposition}

\theoremstyle{definition}
\newtheorem{defn0}[thm0]{Definition}
\newtheorem{nota0}[thm0]{Notation}
\newtheorem{exa0}[thm0]{Example}

\newtheorem{rem0}[thm0]{Remark}
\newtheorem{claim0}[thm0]{Claim}

\newcommand{\elena}[1]{{\color{blue}#1}}


\author[E.~Angelini]{Elena Angelini}
\address{Dipartimento di Ingegneria dell'Informazione e Scienze Matematiche, Universit\`a di Siena, Italy}
\email{elena.angelini@unisi.it}

\author[L.~Chiantini]{Luca Chiantini}
\address{Dipartimento di Ingegneria dell'Informazione e Scienze Matematiche, Universit\`a di Siena, Italy}
\email{luca.chiantini@unisi.it}
\thanks{The authors are members of the Italian GNSAGA-INDAM and are supported by the Italian PRIN 2015 - 
Geometry of Algebraic Varieties (B16J15002000005)}

\begin{document}

\begin{abstract}
We describe a new method to determine the minimality and identifiability of a Waring decomposition $A$ of a specific form (symmetric tensor)
$T$ in three variables. Our method, which is based on the Hilbert function of $A$, can distinguish between forms in the span
of the Veronese image of $A$, which in general contains both identifiable and not identifiable points, depending on the choice
of coefficients in the decomposition. This makes our method applicable for all values of the length $r$ of the decomposition, from $2$ up to the generic rank,
a range which was not achievable before. Though the method in principle can handle all cases of specific ternary forms, 
we introduce and describe it in details for forms of degree $8$.
\end{abstract}

\maketitle

\section{Introduction}
The paper is devoted to the analysis of the identifiability of a Waring  decomposition of a symmetric tensor over $\C$. A symmetric tensor
$T\in S^d\C^{n+1}$ is equivalent to a homogeneous polynomial (form) of degree $d$ in $n+1$ variables, and a Waring decomposition of
$T$ (of length $r$) corresponds to an expression $T=\sum_{i=1}^r L_i^d$, where the $L_i$'s are linear forms. The (Waring) rank of $T$ is the minimal $r$ for
which the decomposition exists, and $T$ is identifiable if the linear forms $L_i$'s appearing in a minimal decomposition are
unique, up to scalar multiplication.

The identifiability of symmetric tensors is relevant for many applications. We refer to the introductions in
 \cite{COttVan17b}, \cite{AnandkumarGeHsuKakadeTelgarsky14},
\cite{AllmanMatiasRhodes09}, \cite{AppellofDavidson81}, and to the many papers cited there,
for an account on how the uniqueness of a decomposition of a tensor $T$  is a fundamental property for algorithms  
in signal processing, image reconstruction, artificial intelligence, statistical mixture models, etc.

In particular, in several concrete cases, one can find a Waring decomposition of a given $T$, either by heuristic computations or by construction.
So the problem is to find criteria which determine whether a given decomposition has minimal cardinality and whether it is unique or not.

The problem was classically solved for binary forms by Sylvester. Thus we mainly focus on the case of ternary forms.

Write $r_d$ for the \emph{generic rank} of ternary forms of degree $d$, i.e. the rank realized outside a 
Zariski closed subset of the space of all degree $d$ forms. By \cite{AlexHir95} and  \cite{COttVan17a}, we know that a \emph{general}
form of rank $r<r_d$ is identifiable, as soon as $d>4$. We will describe below a method to determine, for a given \emph{specific} form $T$,
whether or not a given decomposition of any length $r<r_d$ has minimal cardinality and is unique (up to rescaling).

We will take the projective point of view, to attack the problem. Thus, a linear form $L$ is identified with a point of the projective space
$\Pj^2$ of linear forms. If $v_d:\Pj^2\to \Pj^N$ denotes the Veronese map of degree $d$,  then a {\it decomposition} of $T$ corresponds  to a finite set 
of linear forms $A\subset\Pj^2$ such that $T$ belongs to the linear span of $v_d(A)$.

The most celebrated (and applied) method for detecting the identifiability of a tensor has been introduced by Kruskal \cite{Kruskal77}.
Geometrically, it can be rephrased in terms of the \emph{Kruskal's rank} of the finite set $A$. Several extensions of the Kruskal's 
criterion are available, e.g. the Reshaped Kruskal's Criterion introduced in \cite{COttVan17b}, see Theorem \ref{thm:kr} below. 
Similar analysis  can be found in  papers by Mourrain and Oneto \cite{MourOneto} and Ballico \cite{Ball19}.
Another analysis, based on catalecticant maps and inverse systems, can be found in \cite{MassaMellaStagliano18}.
Yet, all these methods can work, for theoretical reasons, only for values of $r$ which remain far below the generic rank
(see e.g. Proposition 4.9 of \cite{COttVan17a},  for the Reshaped Kruskal's Criterion, as well as the statements of Theorem 1.1 of
\cite{Ball19} and Theorem 2.17 of \cite{MourOneto}).

There is indeed an intrinsic weakness, both in  the Kruskal's and in the catalecticant approaches: they only consider
{\it projective} properties of  the sets $A$ and $v_d(A)$, and not the specific tensor $T$ in the span of $v_d(A)$. 
This means that the methods cannot distinguish between two forms $T_1,T_2\in \langle v_d(A)\rangle$, i.e. 
forms that can be written as $T_1= \sum_{i=1}^r a_iL_i^d$, $T_2=\sum_{i=1}^r b_iL_i^d$ for a (projectively) different
choice of the coefficients. Thus, when the span of $v_d(A)$ contains both tensors for which
$A$ is minimal and unique and tensors for which $A$ is not, then the previous criteria will not apply.
In other words, the previous criteria can determine the identifiability of  $T$ only if \emph{all} the tensors in the span
of $v_d(A)$ (except those spanned by a proper subset) are identifiable. It turns out (see e.g. Example \ref{sharpn}) that
even if $A$ is generic, as soon as the cardinality $r$ approaches the generic value $r_d$  one can find, in the span of
$v_d(A)$, both points for which $A$ is minimal and unique and points for which $A$ is not.
In geometric terms, what happens is that, even for $A$ general, the span of $v_d(A)$ contains  points
in which two folds of the $r$-secant variety cross each other (singular, non-normal points, see Remark \ref{nonsing}).
This implies that the previous analysis \emph{cannot} determine the identifiability of $T$, as soon as $r$ grows.

In section \ref{beyond} we compute the maximal $r_0<r_d$, as a function of $d$, for which an analysis of the decomposition $A$ alone
can determine the identifiability of a ternary form $T$ (see Theorem \ref{range} and Theorem \ref{ranger}  and  Example \ref{sharpn}).
Other analysis (see e.g. the procedure
 described  by Domanov and De Lathauwer in \cite{DomaLath17}) could in principle take into account not only $A$ but also
 the coefficients of the decomposition of $T$,
 but their range of applicability remains, as far as we know, under the bound $r_0$ above, hence far below the generic rank.

So, in order to analyze the minimality and identifiability of a decomposition $A$ of cardinality greater than the bound $r_0$, a deeper analysis
is needed. The analysis must be able to distinguish between different tensors $T_1,T_2$ belonging to the span of the images of the same linear forms
$v_d(L_i)$'s. One of the main targets of the paper is the introduction of a
procedure for such an analysis (see Section 4 below). Even if not directly, implicitly our analysis takes into account also the coefficients of the decomposition
$\sum_{i=1}^r a_iL_i^d$.  Indeed, we explain the procedure in the specific case of ternary
forms of degree $8$, the lowest degree for which our construction becomes effective, see \cite{AngeCMazzon}. We claim, however,
that with the same method one can analyze the behavior of ternary forms of any degree.

The analysis that we propose in Section 4, which is the core of the paper,
 is based on the study of the Hilbert function and a resolution of the ideal of the set $A$. The Hilbert function (see Definition \ref{Hilbdef} below) is a 
central tool for the study of the geometry of finite subsets of projective spaces. It is known that there are connections between properties of the Hilbert
function of $A$ and the identifiability of a tensor $T$ in the span of $v_d(A)$ (see \cite{COttVan17b}, \cite{AngeCVan18}, \cite{AngeCMazzon}).

We are able, by testing the Hilbert function and a resolution of the ideal of $A$, 
to produce algorithms (see \ref{algor} below)  that can guarantee the uniqueness of a given decomposition of a given ternary form $T$, 
in principle for all degrees $d$ and {\it all} values of $r$ smaller than the generic rank $r_d$. As far as we know, this is the first example of 
an analysis which can determine the identifiability of $T$, for all values of $r$ up to the generic rank.

In addition, when a second decomposition $B$ of the same cardinality $r$ exists for $T$,
then our method also indicates how one can construct the second decomposition. Let us stress that Example \ref{sharpn} has the following geometric 
consequence: it shows how our algorithm can detect points $T$ in the span of $v_8(A)$ with two {\it disjoint} decompositions $A$, $B$ of the same length $14$. 
From a geometric point of view, these examples, whose existence is guaranteed by Example \ref{sharpn}, correspond to points $T$ of the secant variety
$\sigma_{14}(v_8(\Pj^2))$ which are {\it non-normal}, see Example \ref{nonsing}.
As far as we know, this is the first example of such singular points in a general secant space to a variety. Their existence indicates that one cannot hope,
except for few initial cases, to solve the identifiability problem for specific forms exclusively with a local analysis.
 
We describe the theoretical basis of the algorithm, in section \ref{otto}, for the case $d=8$ and for rank $r=14$, the biggest value smaller 
than the generic rank $r_8=15$. 
This is the first numerical case in which, for a {\it general} choice of the set $A$ of cardinality $r$, the general form in the span of $v_d(A)$ is identifiable, 
but the span also contains  forms $T$ having another decomposition $B$ of cardinality $r$ (and $B\cap A=\emptyset$).
We notice that, for us, the word \emph{general} has an effective, computable meaning: $A$ is general if some higher Kruskal's ranks of $A$ are general. 
The algorithm is effective, and requires just to control if a certain $12\times 13$ linear system is not solvable. We give examples  of applications
(Example \ref{exok}), and also discuss its computational complexity.

Our analysis can be extended, under the same guidelines, for higher values of $d$.  The (next) case of ternary forms of degree $9$,
which have several geometric peculiarities, will be the topic of a forthcoming paper.  We stress that the method can analyze
even the case of decompositions $A$ whose Kruskal's ranks are not generic. Since for any given value of $r,d,$ and the Kruskal's ranks,
the shape of a resolution of the ideal of $A$ is different, then one needs to adapt the algorithm to the case under analysis.

With the same approach, in principle we could analyze  also the case of forms in $4,5,\dots$ variables. As our knowledge on the Hilbert functions
of finite sets in $\Pj^3, \Pj^4, \dots$ is (by far) less complete than for sets of points in $\Pj^2$, a precise algorithm for the identifiability
of specific forms in many variables is still unavailable. We  observe that, in this way, the theory of tensors can suggests problems 
 in the  geometry of finite projective sets, whose solution could determine relevant 
 theoretical and practical advances, for our knowledge. 
\smallskip

The paper is structured as follows: in section \ref{sec:notation} we introduce main notation and definitions used throughout the paper and 
we recall the symmetric version of Kruskal's criterion. Moreover, some elementary results about the Hilbert function and the 
Cayley-Bacharach property for finite sets are recalled. By means of these tools, in section \ref{beyond} we describe a new method to 
determine the minimality and identifiability of a Waring decomposition of a specific ternary form of sub-generic rank. This analysis allows us 
to go beyond the range of applicability of Kruskal's approach and can be extended in a natural way to the case of a form with an arbitrary 
number of variables. Finally, in section \ref{otto}, we show how the study of the resolution of a decomposition yields a method to determine
 the identifiability of ternary forms, even when it depends on the coefficients of the decomposition. We do that by analyzing specifically 
 the case of ternary forms of degree $ 8 $.

\section{Preliminaries}\label{sec:notation}

\subsection{Notation} \quad

Let $d,n \in \N$. Let $ \C^{n+1} $ be the space of linear forms in $ x_{0}, \ldots, x_{n} $ and $ S^{d} \C^{n+1}$ the space of forms of degree 
$d$ in $ x_{0}, \ldots, x_{n} $ over $ \C $. \\
Let $ T \in S^{d} \C^{n+1} $. $ T $ is associated to an element of $ \Pj(S^{d} \C^{n+1}) \cong \Pj^{N} $ $( N = \binom{n+d}{d} - 1) $, which, 
by abuse of notation, we denote by $T$. \\
Let $ v_{d}: \Pj^{n} \rightarrow \Pj^{N} $ be the \emph{Veronese embedding} of $ \Pj^{n} $ of degree $ d $, which is given by 
$$ v_{d}([a_{0}x_{0}+ \ldots + a_{n}x_{n}]) = [(a_{0}x_{0}+ \ldots + a_{n}x_{n})^{d}]. $$ 
Let $ A = \{P_{1}, \ldots, P_{\ell(A)}\} \subset \Pj^{n} $ be a finite set of cardinality $ \ell(A) $. We define $v_{d}(A) = 
\{v_{d}(P_{1}), \ldots , v_{d}(P_{\ell(A)})\}$ and we denote by $ \langle v_{d}(A) \rangle $ the linear space spanned 
by $ v_{d}(P_{1}), \ldots, v_{d}(P_{\ell(A)}) $. 

\smallskip

With the above notations we give the following definitions.

\begin{defn0}
Let $ A \subset \Pj^n $ be a finite set. $ A $ \emph{computes} $ T $ if $ T\in \langle v_{d}(A) \rangle$, 
the linear space spanned by the points of $ v_{d}(A) $.
\end{defn0}

\begin{defn0}
Let $ A \subset \Pj^n $ be a finite set which computes $ T $. $ A $ is \emph{non-redundant} if we cannot find a proper subset $ A' $ of $ A $ 
such that $T \in \langle v_{d}(A')\rangle $.
\end{defn0}

\begin{rem0}\label{rem:indep}
If $ A \subset \Pj^n $ is a finite set that computes $T$ and it is non-redundant, then the points of $ v_{d}(A) $ are 
linearly independent, i.e.,
$$ \dim(\langle v_{d}(A)\rangle) = \ell(A) -1.  $$
\end{rem0}

Moreover we introduce the following:

\begin{defn0}
The \emph{rank} of $ T $ is $ r = \min \,\{\ell(A) \, |\, T \in \langle v_{d}(A)\rangle\} $. A finite set $ A \subset \mathbb{P}^{n} $
 \emph{computes the rank} of $ T $ if $ A $ computes $ T $, it is non-redundant and $ \ell(A) = r $.
\end{defn0}

\begin{defn0}
$ T $ of rank $ r $ is \emph{identifiable} if there exists a unique $ A $ computing the rank of $ T $. 
\end{defn0}

\subsection{Kruskal's criterion for symmetric tensors}\label{sec:Kr}

\begin{defn0}
The \emph{d-th Kruskal's rank} of a finite set $ A \subset \mathbb{P}^{n} $ is 
$$ k_{d}(A) = \max \,\{k \, | \, \forall \, A' \subset A, \, \ell(A') \leq k, \, \dim \langle v_{d}(A')\rangle = \ell(v_{d}(A')) - 1 \}. $$
\end{defn0} 

\begin{rem0}\label{maxKrank}
For any $ d $, it holds that $ k_{d}(A) \leq \min\{N+1,\ell(A)\} $. Moreover, 
if $k_d(A)=\min \min\{N+1,\ell(A)\}$ is maximal, then for all $A'\subset A$ the Kruskal's rank
$k_d(A')$ is also maximal. 

If $ A $ is sufficiently general, then $ k_{d}(A) = \min\{N+1,\ell(A)\} $ (see e.g. Lemma 4.4 of \cite{COttVan17b}.
\end{rem0}

The Kruskal's rank is fundamental in the statement of the reshaped Kruskal's criterion.
 
\begin{thm0}[Reshaped Kruskal's Criterion, see \cite{COttVan17b}]\label{thm:kr}
Let $ T \in \Pj(S^{d} \C^{n+1}) $ with $ d \geq 3 $ and let $ A \subset \mathbb{P}^{n} $ be a non-redundant set computing $ T $.  
Assume that $ d = d_{1}+d_{2}+d_{3} $ with $ d_{1} \geq d_{2} \geq d_{3} \geq 1 $. If 
\begin{equation}\label{eq:Kr}
\ell(A) \leq \frac{k_{d_{1}}(A)+k_{d_{2}}(A)+k_{d_{3}}(A)-2}{2}
\end{equation}
then $ T $ has rank $ \ell(A) $ and it is identifiable.
\end{thm0}

\subsection{The Hilbert function for finite sets in $ \Pj^{n} $}

\begin{defn0}
The \emph{evaluation map} of degree $ d $ on a ordered finite set of vectors $ Y = \{Y_{1}, \ldots, Y_{\ell}\} \subset \mathbb{C}^{n+1} $ is the linear map given by
$$ ev_{Y}(d): Sym^{d}\mathbb{C}^{n+1} \longrightarrow \mathbb{C}^{\ell} $$
$$ ev_{Y}(d)(F) = (F(Y_{1}), \ldots, F(Y_{\ell})). $$
\end{defn0} 

\begin{defn0}\label{Hilbdef}
Let $ Y $ be a set of homogeneous coordinates for a finite set $ Z $ of $ \mathbb{P}^{n} $. The \emph{Hilbert function} of $ Z $ is the map
$$ h_{Z}: \mathbb{Z} \longrightarrow \mathbb{N} $$
such that $ h_{Z}(j) = 0$, for $j < 0$, $ h_{Z}(j) = rank (ev_{Y}(j))$, for $ j \geq 0. $
\end{defn0} 

\begin{rem0}\label{spans} Take the notation of the previous definition. Since elements of the kernel of the evaluation map $ev_Y(1)$ 
correspond to the equations of hyperplanes vanishing at $Y$, it turns out that
$h_Z(1)$  is the (affine) dimension of the linear space spanned by $Z$.

Since elements of the kernel of the evaluation map $ev_Y(d)$ correspond to the equations of hypersurfaces of degree $d$ vanishing at $Y$, which
in turn correspond to the equations of hyperplanes vanishing at $v_d(Z)$, thus it corresponds to the (affine) dimension of the span $\langle v_d(Z)\rangle$.
\end{rem0}

\begin{defn0}
The \emph{first difference of the Hilbert function} $ Dh_{Z} $ of $ Z $ is 
$$ Dh_{Z}(j) = h_{Z}(j) - h_{Z}(j-1), \, j \in \mathbb{Z}. $$
\end{defn0}

\begin{rem0}\label{HilbElem} We collect here some useful elementary properties of $ h_{Z} $ and $ Dh_{Z} $:
\begin{itemize}
\item[$1.$] $ Dh_{Z}(j) = 0$, for $ j < 0$;
\item[$2.$] $ h_{Z}(0) = Dh_{Z}(0) = 1$;
\item[$3.$] $ Dh_{Z}(j) \geq 0 $, for all $ j $;
\item[$4.$] $ h_{Z}(i) = \sum_{0\leq j \leq i} Dh_{Z}(j) $; 
\item[$5.$] $ h_{Z}(j) = \ell(Z) $, for all $ j \gg 0$; 
\item[$6.$] $ Dh_{Z}(j) = 0$, for $ j \gg 0 $;
\item[$7.$] $ \sum_{j} Dh_{Z}(j) = \ell(Z) $. 
\end{itemize}

The proofs appear, sparsely, in several sources (\cite{IK}, \cite{GerMaroRoberts83}, \cite{Migliore}).
They are collected in Lemma 6.1 of \cite{C19}).
\end{rem0}

\begin{prop0}\label{prop:inclusion}
If $ Z' \subset Z $, then, for any $ j \in Z $, it holds that 
$$ h_{Z'}(j) \leq h_{Z}(j), \, Dh_{Z'}(j) \leq Dh_{Z}(j). $$  
\end{prop0}
\begin{proof} Well known. See e.g. Proposition 6.3 of \cite{C19}.
\end{proof}

\begin{prop0}\label{prop:nonincr} 
If there exists $  i >0 $ such that $ Dh_{Z}(i) \leq i $, then 
$$ Dh_{Z}(i) \geq Dh_{Z}(i+1). $$ 
Therefore, if $ Dh_{Z}(i) = 0 $, then $ Dh_{Z}(j) = 0 $ for any $ j \geq i $.
\end{prop0}
\begin{proof} Well known. See e.g. Proposition 6.4 of \cite{C19}.
\end{proof}

The following theorem, which provides one of the main tools for our analysis, was proved by E. Davis for
points in $\Pj^2$ (\cite{Davis85}), and then extended to points in any projective space in Theorem 3.6 of \cite{BigaGerMig94}.

\begin{thm0}[Davis 1985]\label{thm:Davis}
Let $ Z \subset \mathbb{P}^{2} $ be a finite set. Assume that:
\begin{itemize}
\item[$1.$] $ Dh_{Z}(j) = j+1 $ for $ j \in \{0, \ldots, i-1\} $ and $ Dh_{Z}(i) \leq i $;
\item[$2.$] $ Dh_{Z}(j_{0}) = Dh_{Z}(j_{0}+1) = e $ for some $ j_{0} \geq i-1 $.
\end{itemize}
Then $ Z = Z_{1} \cup Z_{2} $, where $ Z_{1} $ lies on a curve of degree $e$ of $\mathbb{P}^{2} $, $Dh_Z(j)=Dh_{Z_1}(j)$
for all $j\geq j_0$, and, for any $  j \in \{0, \ldots, j_{0} - e - 1\} $, 
$ Dh_{Z_{2}}(j) = Dh_{Z}(e+j) - e. $
\end{thm0}

\begin{nota0} Let $ Z \subset \mathbb{P}^{2} $ be a finite set and let $d \in \N$. We pose
$$h^1_Z(d)= \ell(Z)-h_Z(d) = \sum_{j=d+1}^\infty Dh_Z(j).$$
\end{nota0}

We recall the following results:
\begin{prop0}\label{d+1}  (See Lemma 1 of \cite{BallBern12a}).
Let $ T \in S^{d}\C^{n+1} $ and let $ A, B \subset \Pj^{n} $ be non-redundant finite sets computing $ T $. 
Put $ Z = A \cup B \subset \Pj^{n}$. Then $ Dh_{Z}(d+1) > 0 $.
\end{prop0}

\begin{prop0}\label{cap} (See formula (15) of \cite{AngeCVan18}, which 
is stated only for degree $d=4$, but whose proof works indeed for any $d$).
Let $A,B\subset \Pj^n$ be finite sets and set $Z=A\cup B$.
For any $d \in \N$, 
$$\dim(\langle v_d(A)\rangle\cap \langle v_d(B)\rangle) = \ell(A \cap B) - 1 + h^1_Z(d).$$
\end{prop0}

As a consequence of Theorem \ref{thm:Davis} and Proposition \ref{cap}, we get the following:
\begin{prop0}\label{Dav}
Let $ T \in S^{d}\C^{3} $ and let $ A \subset \Pj^{2} $ be a non-redundant finite set computing $ T $. 
Then, there is no other $ B \subset \Pj^{2} $ non-redundant finite set computing $ T $ with $ A \cap B = \emptyset $, 
$ \ell(B) \leq \ell(A) $ and such that, if $ Z = A \cup B \subset \Pj^{2}$, then:
\begin{itemize}
\item[$1.$] $ Dh_{Z}(j) = j+1 $ for $ j \in \{0, \ldots, i-1\} $ and $ Dh_{Z}(i) \leq i $;
\item[$2.$] $ Dh_{Z}(j_{0}) = Dh_{Z}(j_{0}+1) = e < i $ for some $ j_{0} > i-1 $.
\end{itemize}
\end{prop0}
\begin{proof}
Assume that such $ B $ exists. Then, by Theorem \ref{thm:Davis}, there exists a proper subset $ Z' $ of $ Z $ contained in a plane curve of degree $ e $, and satisfying $ h^{1}_{Z}(d) = h^{1}_{Z'}(d) $. By Proposition \ref{cap} $ h^{1}_{Z}(d) > 0 $, being $ A $ and $ B $ non-redundant decompositions for $ T $. 
Set $ Z' = A' \cup B' $, with $ A' \subset A, B' \subset B$. We have that $ A' \cap B' = \emptyset $ and that 
$ A' \subsetneq A $ or $ B' \subsetneq B $. Therefore, by Proposition \ref{cap},
$$\dim(\langle v_d(A')\rangle\cap \langle v_d(B')\rangle) = - 1 + h^1_{Z'}(d) = \dim(\langle v_d(A)\rangle\cap \langle v_d(B)\rangle)  $$
and so $ T \in \langle v_d(A')\rangle\cap \langle v_d(B')\rangle $, which violates the non-redundantity assumption on $ A $ and $ B $, 
depending on whether $ A' \subsetneq A $ or $ B' \subsetneq B$.
\end{proof}

\subsection{The Cayley-Bacharach property for finite sets in $ \mathbb{P}^{n} $} \label{sec:CB}

\begin{defn0}
A finite set $Z\subset \Pj^{n}$ satisfies the \emph{Cayley-Bacharach property in degree $d,$} $\mathit{CB}(d)$, if, for all $ P \in Z$, 
it holds that every form of degree $d$ vanishing at $ Z\setminus\{ P\}$ also vanishes at $P$.
\end{defn0} 

\begin{exa0} \quad
\begin{itemize}
\item[$1.$] Let $ Z \subset \mathbb{P}^{2} $ be a set of $6$ general points. Then 
$$ \begin{tabular}{c|ccccc}
$j$ & $0$ & $1$ & $2$ & $ 3 $ & $ \dots $ \\  \hline
$h_Z(j)$ & $1$ & $3$ &   $6$ & $ 6 $ & $ \dots $ \cr
$Dh_Z(j)$ & $1$ & $2$ &   $3$ & $ 0 $ & $ \dots $ \cr
\end{tabular} $$
and $ Z $ has $\mathit{CB}(1)$ but not $\mathit{CB}(2)$. 
\item[$2.$] Let $ Z \subset \mathbb{P}^{2}$ be a set of $6$ points on an irreducible conic. Then 
$$ \begin{tabular}{c|cccccc}
$j$ & $0$ & $1$ & $2$ & $ 3 $ & $4 $ & $ \dots $ \\  \hline
$h_Z(j)$ & $1$ & $3$ &   $5$ & $ 6 $ & $ 6 $ & $ \dots $ \cr
$Dh_Z(j)$ & $1$ & $2$ &   $2$ & $ 1 $ & $ 0 $ & $ \dots $ \cr
\end{tabular} $$
and $ Z $ has $\mathit{CB}(2)$ and $\mathit{CB}(1)$. 
\item[$3.$] Let $ Z \subset \mathbb{P}^{2} $ be a set of $6$ points, of which $ 5 $ aligned. Then
$$ \begin{tabular}{c|ccccccc}
$j$ & $0$ & $1$ & $2$ & $ 3 $ & $4 $ & $ 5 $ & $ \dots $ \\  \hline
$h_Z(j)$ & $1$ & $3$ &   $4$ & $ 5 $ & $ 6 $ &  $ 6 $ & $ \dots $ \cr
$Dh_Z(j)$ & $1$ & $2$ &   $1$ & $ 1 $ & $ 1 $ & $ 0 $ & $ \dots $ \cr
\end{tabular} $$
and $ Z $ has not $\mathit{CB}(1)$.
\end{itemize}
\end{exa0}

Some fundamental consequences of the Cayley-Bacharach property are listed below. 

\begin{prop0} \label{CBh1} If $Z$ satsfies the property $\mathit{CB}(d)$, then for any proper subset $Z'\subset Z$
we have $h^1_{Z'}(d)<h^1_{Z}(d)$.
\end{prop0}
\begin{proof} Since $\sum_{i=1}^\infty Dh_Z(i) =\ell(Z) > \ell(Z')=\sum_{i=1}^\infty Dh_{Z'}(i)$, if $h^1_{Z'}(d)=h^1_Z(d)$
then $\sum_{i=1}^d Dh_Z(i) >\sum_{i=1}^d Dh_{Z'}(i)$, i.e. $h_Z(d)>h_{Z'}(d)$. Thus, in degree $d$, the homogeneous ideal of $Z'$ has dimension
bigger than the homogeneous ideal of $Z$, which means that there exists some form of degree $d$ containing $Z'$ and
not containing $Z$.
\end{proof}

\begin{thm0}[Angelini, Chiantini, Vannieuwenhoven 2018, \cite{AngeCVan18}]\label{thm:extGKR}
If $ Z $ has $\mathit{CB}(d)$, then, for any $ j \in \{0, \ldots, d+1\} $, it holds that
\begin{equation}\label{eq:CB}
Dh_{Z}(0) + \ldots + Dh_{Z}(j) \leq Dh_{Z}(d+1-j) + \ldots + Dh_{Z}(d+1).
\end{equation}
\end{thm0}

As in \cite{AngeCVan18}, the Cayley-Bacharach property is relevant in our analysis since it holds for sets $Z=A\cup B$,
where $A,B$ are two different \emph{non-redundant, disjoint} decompositions of a form $T$. 

Next proposition is essentially contained in \cite{AngeCMazzon} (Lemma 5.3). 

\begin{prop0}\label{CBconseq}
Let $ T \in S^{d}\C^{n+1} $ and let $ A \subset \Pj^{n} $ be a non-redundant finite set computing $ T $. 
Let $ B \subset \Pj^{n} $ be another non-redundant finite set computing $ T $
and assume $ A \cap B = \emptyset $. Then $Z=A\cup B$ satisfies the Cayley-Bacharach property $\mathit{CB}(d)$. 
\end{prop0}
\begin{proof}
Assume that $Z$ does not satisfy $\mathit{CB}(d)$. Then there exists $P\in Z=A\cup B$ such that the ideal of $Z\setminus\{P\}$ is strictly bigger than the ideal
of $Z$ in degree $d$. This implies that:
$$h_Z(d)=\sum_{i=0}^{d}Dh_Z(i) >  \sum_{i=0}^{d}Dh_{Z\setminus\{P\}}(i)=h_{Z\setminus\{P\}}(d).$$
Since $Dh_Z(i)\geq Dh_{Z\setminus\{P\}}(i)$ for all $i$ (Proposition \ref{prop:inclusion}) and: 
$$\ell(Z)=\sum_{i=0}^{\infty}Dh_Z(i) =1+\sum_{i=0}^{\infty}Dh_{Z\setminus\{P\}}(i)=
1+\ell(Z\setminus\{P\}), $$
then necessarily $h^1_Z(d)=h^1_{Z\setminus\{P\}}(d)$, so that, by Proposition \ref{cap}:
{\small $$\dim(\langle v_d(A)\rangle\cap \langle v_d(B)\rangle) = h^1_Z(d)-1 = h^1_{Z\setminus\{P\}}(d)-1=
\dim(\langle v_d(A\setminus\{P\})\rangle\cap \langle v_d(B\setminus\{P\})\rangle).$$}
Thus $T\in \langle v_d(A\setminus\{P\})\rangle\cap \langle v_d(B\setminus\{P\})\rangle$, which contradicts 
 the assumption that both $A$ and $B$ are non-redundant.
\end{proof}

\section{Beyond the Kruskal's bound for forms in three variables}\label{beyond}

In this section we prove a sharp criterion  which determines the identifiability of a form $T$ of degree $d$ in $3$ variables,
in terms of linear algebraic invariants on the coordinates of the points of a decomposition of $T$.
\smallskip

Following the general notation, let  $A\subset\Pj^2$ be a non-redundant set which computes  $T$. Put $r=\ell(A)$.  We want to find
a criterion, based on the geometric properties of $A$, which guarantees that $T$ is identifiable of rank $r$. The criterion should be
effective on an ample collection of decompositions.

\begin{thm0}\label{range} The form $T$  is identifiable of rank $r$  if one of the following holds:

\begin{itemize}
\item $d=2m$ is even, $k_{m-1}(A)= \min\{\binom{m+1}2,r\}$,  $h_A(m)=r\leq \binom{m+2}2-2$;

\item $d=2m+1$ is odd, $k_{m}(A)=\min\{\binom{m+2}2,r\}$, $h_A(m+1)=r\leq \binom{m+2}2 +\lfloor {m \over 2} \rfloor $.
\end{itemize}
\end{thm0}

The numerical assumptions on $k_{m-1}(A)$, $h_A(m)$, $k_{2e}(A)$, $k_{2e+1}(A)$, $h_A(2e+1)$, $h_A(2e+2)$ imply that these
values are maximal, for a set $A$ of $r$ points.
Thus, the assumptions of Theorem \ref{range} are expected to hold, provided that $A$ is a sufficiently general set of points (Remark \ref{maxKrank}).

\begin{rem0} 
The complexity of of the algorithm for computing the Kruskal's ranks in the assumptions of Theorem \ref{range}
can be computed as follows.

When $d=2m$, if one puts as rows of a matrix $M_m$ (resp. $M_{m-1}$) a set of homogeneous coordinates of the points of $v_m(A)$
(resp. $v_{m-1}(A)$), then $h_A(m)=r$ simply means that the matrix $M_m$ has full rank $r$. This in general needs the computation
of one $r\times r$ minor. In order to control that $k_{m-1}(A)=\min\{\binom{m+1}2,r\}$,
one has to compute (in general)  just one $r\times r$ determinant, when  $\binom{m+1}2\geq r$.
On the other hand, when $\binom{m+1}2<r$, the computation of $k_{m-1}(A)$ requires the computation
of all the minors of $M_{m-1}$ obtained by taking any subset of $\binom{m+1}2$ rows. Since $M_{m-1}$ has $r$ rows
and $r\leq \binom{m+1}2+m$,  then one must compute the maximality of the rank of (at worst) $r^m/m!$ matrices of type $ \binom{m+1}2\times  \binom{m+1}2$.
\end{rem0}
\medskip 

Before proving Theorem \ref{range}, we want to point out the following remark, which will be useful in several arguments, to handle the case of tensors 
with two decompositions $ A,B $ with $ A \cap B \not= \emptyset $:

\begin{rem0}\label{redint}
Assume that $ T $ is a form of degree $ d $ in $ n+1 $ variables, computed by a non-redundant finite set $ A = \{P_{1}, \ldots, P_{r}\}\subset \Pj^{n} $.\\
\noindent Let $B$ be another non-redundant decomposition of $T$ with $ s = \ell(B)\leq r$ and define $Z=A\cup B$. \\
\noindent If the intersection $A\cap B$ is not empty, then we can reorder the points of $A$ so that $B=\{P_1,\dots,P_j,P'_{j+1},\dots,P'_s\}$ with $j\geq 1$ 
and $P'_i\notin A$ for $i=j+1,\dots,s$. Then for any choice of representatives (i.e. coordinates) $T_1,\dots T_r, T'_{j+1},\dots,T'_{s}$ for the projective points 
$v_d(P_1),\dots,v_d(P_r),v_d(P'_{j+1}),\dots, v_d(P'_s)$ respectively, there are non-zero scalars $a_i$'s, $b_i$'s such that
$$\begin{matrix}
T & = &a_1T_1+\dots+a_rT_r \\ T & = & b_1T_1+\dots+b_jT_j+b_{j+1}T'_{j+1}+\dots+b_sT'_s.
\end{matrix}$$
Define:
\begin{multline*} T_0=(a_1-b_1)T_1+\dots+(a_j-b_j)T_j+a_{j+1}T_{j+1}+\dots+a_rT_r
\\  = b_{j+1}T'_{j+1}+\dots+b_sT'_s.\end{multline*}
Now $T_0$ has the two decompositions $A$ and $B'=\{P'_{j+1},\dots, P'_s\}$, which are disjoint.  If $B'$ is redundant, then after 
rearranging the points, we may assume
$T_0=c_{j+1}T'_{j+1}+\dots+c_tT'_t$ for some  $t<s$, so that:
\begin{multline*} T=b_1T_1+\dots +b_jT_j+T_0=  b_1T_1+\dots +b_jT_j+c_{j+1}T'_{j+1}+\dots+c_tT'_t,
\end{multline*}
against the  fact that $B$ is non-redundant. Thus $B'$ must be non-redundant. \\
If $ A $ is non-redundant, define $ A' = A $. \\
If $A$ is redundant, since the points $v_d(P_1),\dots ,v_d(P_r)$ are linearly independent (Remark \ref{rem:indep}), then some coefficient $(a_i-b_i)$ is $0$. 
In this case, we may assume $(a_i-b_i)=0$ if and only if $i=1,\dots,q\leq j$, so we get a non-redundant decomposition $A'=\{P_{q+1},\dots,P_r\} $ of $ T_{0} $. \\
In conclusion, we find that $ T_{0} $ has two different decompositions non-redundant $A', B'$, with $ A' \subset A $ and $ \ell(B') \leq \ell(A') $, $ \ell(B') < r $.
\end{rem0}

{\it Proof of Theorem \ref{range}}
We prove the statement by induction on $ r$, the case $ r = 1 $ being trivial. Thus we may assume $ d \geq 4 $. \\
\indent With the notation of Remark \ref{redint}, the existence of another non-redundant decomposition $ B $ of $T$ with $\ell(B)\leq r$ and intersecting $ A $, 
implies the existence of a new form $ T_{0} $ admitting two disjoint non-redudant decompositions $ A', B' \subset \Pj^{2} $, with $ A' \subset A $. \\
\indent Since $A'\subset A$, then the evaluation map in degree $m$ (respectively $ m+1 $, depending on the parity of $ d $) surjects for $A'$ (see Remark \ref{maxKrank}). 
 Moreover $k_{m-1}(A')= \min\{\binom{m+1}2,r'\}$  (respectively $k_{m}(A')= \min\{\binom{m+2}2,r'\}$, depending on the parity of $ d $), by Remark \ref{maxKrank}. \\
\indent Assume $ \ell(A') < \ell(A) $, since $T_0$ has a second non-redundant decomposition $B'$  of length $\ell(B')\leq \ell(A')$, by induction we get a contradiction. \\
\indent Thus $A = A'$ and $B'$ are two non-redundant decompositions of $T_0$, with $\ell(B')<r$.
 By replacing $T,B$ with $T_0,B'$ respectively, \emph{we can thus reduce ourselves to prove the claim
 only in the case $A\cap B=\emptyset.$} 
In order to do that, we need to distinguish two cases. \\
\indent First, assume $d=2m$. Since $A\cap B=\emptyset$, then by Proposition \ref{CBconseq} $Z$ satisfies $\mathit{CB}(d)$. It follows by Proposition 
\ref{prop:inclusion} and by Theorem \ref{thm:extGKR}:
$$ \sum_{i=0}^{m-1} Dh_A(i)\leq \sum_{i=0}^{m-1} Dh_Z(i)\leq \sum_{i=m+2}^{d+1} Dh_Z(i).$$
Moreover, our hypotheses imply that 
\begin{multline*} r- \sum_{i=0}^{m-1} Dh_A(i) = r - h_A(m-1) = r -  k_{m-1}(A) = 
              \left\{
                \begin{array}{ll}
                  0 & if \, r \leq \binom{m+1}2\\
                  r - \binom{m+1}2 & if \, r > \binom{m+1}2 \\
                \end{array}
              \right.  \leq \\  
                \left\{
                \begin{array}{ll}
                  0 & if \, r \leq \binom{m+1}2\\
                  \binom{m+2}2 - 2 - \binom{m+1}2 & if \, r > \binom{m+1}2 \\
                \end{array}
              \right.
               \leq  m-1 < m.
\end{multline*} 
from the previous formulas we get that:
$$Dh_Z(m)+Dh_Z(m+1)\leq 2r-\sum_{i=0}^{m-1} Dh_Z(i)-\sum_{i=m+2}^{d+1} Dh_Z(i)\leq 
 2r-2\sum_{i=0}^{m-1} Dh_A(i)< 2m,$$
i.e. either $Dh_Z(m+1)< m$, or $Dh_Z(m)< m$ in which case, by Proposition \ref{prop:nonincr}, we conclude again that  $Dh_Z(m+1)\leq Dh_Z(m)< m$.
It follows from Proposition \ref{prop:nonincr} that  $Dh_Z(i)\geq Dh_Z(i+1)$ for $i\geq m+1$. If $Dh_Z(i)>Dh_Z(i+1)$ for $i=m+1\dots,d$, until it reaches $0$, 
then  we get $Dh_Z(d+1)= 0$, a contradiction.\\
\indent Assume now that $d=2m+1$. Just as above, one can show that the following inequalities hold:
 $$ \sum_{i=0}^{m} Dh_A(i)\leq \sum_{i=0}^{m} Dh_Z(i)\leq \sum_{i=m+2}^{d+1} Dh_Z(i)$$
$$ r- \sum_{i=0}^{m} Dh_A(i) = 
              \left\{
                \begin{array}{ll}
                  0 & if \, r \leq \binom{m+2}2\\
                  r - \binom{m+2}2 & if \, r > \binom{m+2}2 \\
                \end{array}
              \right.  \leq 
                \left\{
                \begin{array}{ll}
                  0 & if \, r \leq \binom{m+1}2\\
                  \lfloor {m \over 2} \rfloor & if \, r > \binom{m+1}2 \\
                \end{array}
              \right.
               \leq  \lfloor {m \over 2} \rfloor.
$$  
\noindent Therefore we get that:
\begin{multline*} Dh_Z(m+1)\leq 2r-\sum_{i=0}^{m} Dh_Z(i)-\sum_{i=m+2}^{d+1} Dh_Z(i)\leq \\
 2r-2\sum_{i=0}^{m} Dh_A(i)\leq 2 {\lfloor {m \over 2} \rfloor}< m+1.\end{multline*}
It follows from proposition \ref{prop:nonincr} that  $Dh_Z(i)\geq Dh_Z(i+1)$ for $i\geq m+1$.
If $Dh_Z(i)<Dh_Z(i-1)$ for $i=m+1,\dots,d+1$ until it reaches $0$, then we get $Dh_Z(d+1)= 0$, a contradiction.
Thus there exists $j\geq m+1$, $j\leq d$, such that $Dh_Z(j)=Dh_Z(j+1)\leq Dh_Z(m+1) < m+1$. By Proposition \ref{Dav} we get the contradiction.
\hfill\qed
\medskip

\begin{exa0}\label{sharpn} We prove that the previous bounds are sharp.

Assume that $d=2m$. Take a general set $A$ of $r = \binom{m+2}2-1$ points in $\Pj^2$. The generality of $A$ implies that
$Dh_A(i)=i+1$ for $i=0,\dots m-1$, $Dh_A(m)=m=Dh_A(m-1)$, so that $A$ is contained in a curve $C$ of degree $m$; 
moreover $A$ is in {\it uniform position} (i.e. the Hilbert functions of two subsets of $A$ of the same cardinality are equal), 
so that $C$ is irreducible; finally the ideal of $A$ is generated in degree $m+1$
(all these properties can be found in \cite{GerMaro84} and \cite{GerMaroRoberts83}). It follows by Proposition 4.1 of \cite{PeskineSzpiro74} 
that one can find another curve $C'$ of degree $m+3$
containing $A$, such the complete intersection $Z=C\cap C'$ is formed by $m(m+3)=2r$ distinct points. Take $B=Z\setminus A$, so that
also $B$ is a set of $r$ points, disjoint from $A$, and $Z=A\cup B$. By  \cite{Davis85}, we have $Dh_Z(d+1)=1$, $Dh_Z(d+2)=0$, moreover 
 the Cayley -Bacharach property $\mathit{CB}(d)$ holds for $Z$. It follows by Proposition \ref{cap} that 
 $\langle v_d(A)\rangle$ and $\langle v_d(B)\rangle$ meet in one point $T$, which thus has two decompositions
 of length $r$: $A$ and $B$. We can prove that $A$ is non-redundant as follows: assume that $T\in  \langle v_d(A')\rangle$ 
 for some proper subset $A'\subset A$. Then we have a proper subset $Z'=A'\cup B\subset Z$ such that $h^1_{Z}(d)=1=h^1_{Z'}(d)$.
 This contradicts Proposition \ref{CBh1}. (Notice that also $B$ is non-redundant, for a general choice of $A$, $C'$. Indeed the
 situation between $A$ and $B$ is essentially symmetric).

When $d=4e+1$, we get an example of a form of degree $d$ with two non-redundant decompositions of length  
$r= \binom{2e+2}2 +e+1$ by taking a general set of $r$ points and embedding it in a general complete intersection of type $2e+2,2e+2$.

When $d=4e+3$, we get an example of a form of degree $d$ with two non-redundant decompositions of length  
$r= \binom{2e+3}2 +e+1$  by taking a general set of $r$ points and embedding it in a general complete intersection of type $2e+2,2e+4$.
\end{exa0}

The first case in which the previous examples produce a new phenomenon is $d=8$.  General ternary forms of degree $8$ have rank $15$. 
Thus, by \cite{COttVan17a}, the general ternary form of degree $8$ and rank $14$ is identifiable. Yet, for  a general choice of a set $A$ of $14$
points in $\Pj^2$, the span $v_8(A)$ contains (special) points for which the decomposition $A$ is non-redundant, but there exists another decomposition $B$
of length $14$.

We will analyze in details the  identifiability of ternary forms of degree $8$ in section \ref{otto}.

\begin{exa0} In the statement of Theorem \ref{range}, when $d$ is even, i.e. $d=2m$, 
and $r\leq \binom{m+1}2$, then the numerical assumptions hold exactly when $h_A(m-1)=r$. So, there is no need to compute the Kruskal's ranks,
in this case.

On the other hand, when $r$ is big, we cannot drop the assumption $k_{m-1}(A)= \min\{\binom{m+1}2,r\}$,
or substitute it with an assumption on some value of $h_A$.
 
 Namely, take $d=8$, i.e. $m=4$. Fix a general plane cubic curve $\Gamma$ and a general set of $12$ points $P_1,\dots,P_{12}$
 on $\Gamma$. If $P$ is a general point of  $\Pj^2$, the set $A=\{P_1,\dots,P_{12},P\}$ satisfies $h_A(4)=13$, $h_A(1)=3$
 (it satisfies also $h_A(3)=10=\binom{3+2}2$). Notice that $k_3(A)=9<\min\{13,10\}$. We prove that a general form $T$ in the span of $v_8(A)$ is not identifiable.
 
 Indeed assume $T=\sum_{i=1}^{12} a_iv_8(P_i) + av_8(P)$ and set $T'=\sum_{i=1}^{12} a_iv_8(P_i) $. $T'$ is a tensor
 whose (non-redundant)  decomposition $\{P_1,\dots,P_{12}\}$ lies in $\Gamma$. Since $v_8(\Gamma)$ is an elliptic normal curve,
  It is well known (see \cite{CCi06} or \cite{AngeBocciC18}) that $T'$ has a second decomposition $B'\subset\Gamma$ of length $12$.
  Thus $T$ has a second decomposition $B\cup\{P\}$ of length $13$.
  
  Similar examples prove that one cannot relax the assumption on $k_{2e}(A)$ (resp. $k_{2e+1}(A)$) when $d=4e+1$
 (resp. $d=4e+3$), and $r$ is big.
\end{exa0}

One should compare the statement of Theorem \ref{range} with Theorem 2.17 of \cite{MourOneto}, where the authors prove 
that $T$  is identifiable when $d\geq 2\delta(A)+1$, where $\delta(A)$ is the Castelnuovo-Mumford regularity of $A$. 
The Castelnuovo-Mumford regularity of $A$ is the minimum $i>0$ such that $h_A(i)=\ell(A)$, in other words it is the minimum $i>0$
such that $Dh_A(i+1)=0$.
In our case, when $r$ is maximal, the assumptions of Theorem \ref{range} imply that the Castelnuovo-Mumford regularity $\delta(A)$ is $m$ if $d=2m$, 
it is $2e+1$ if $d=4e+1$ and it is $2e+3$ if $d=4e+3$. 
 Thus, Theorem 2.17 of \cite{MourOneto} does not apply, because e.g. in the even case $d=2m<2\delta(A)+1$. 
 From this point of view, Theorem \ref{range} goes beyond the Mourrain-Oneto's result, for the case of three variables.
Notice indeed that, e.g. in the case $d=2m$, the regularity of $A$ in degree $m-1$ implies that $r\leq \binom{m+1}2$,
 so it  is equivalent to  the conditions  $k_{m-1}(A)=r$ and $h_A(m)=r$.
 \smallskip
 
 Notice that Theorem \ref{range} implies in particular that, under the assumptions of the statement, $T$ has rank $r=\ell(A)$.
Indeed, if one is only interested in the fact that $A$ computes the rank of $T$, and not in the uniqueness of $A$, then the statement can be refined.

\begin{thm0}\label{ranger} The decomposition $A$ of $T$ computes the rank of  $T$  if one of the following holds:

\begin{itemize}
\item $d=2m$ is even and $h_A(m)=r (\leq \binom{m+2}2)$:
\item $d=2m+1$, $k_{m}(A)=\min\{\binom{m+2}2,r\}$, $h_A(m+1)=r\leq \binom{m+2}2 +\lceil {m \over 2} \rceil$.
\end{itemize}
\end{thm0}
\begin{proof} The proof is rather similar to the proof of Theorem \ref{range}. We want to exclude the existence of 
another non-redundant decomposition $B$, with $\ell(B)<\ell(A).$ Here $Z=A\cup B$ has cardinality $\ell(Z)<2r$.

By applying Remark \ref{redint} with $ s = \ell(B) < r $ and by arguing as in the proof of Theorem \ref{range}, one can reduce the proof to the case 
$A\cap B =\emptyset$.

Assume $d=2m$. 
Since $r=\ell(A)=h_A(m)$, then $\sum_{i=0}^m Dh_A(i)=r$, hence $\sum_{i=0}^m Dh_Z(i)\geq r$. Then, by Proposition \ref{CBconseq}:
$$ \ell(Z) \geq 2\sum_{i=0}^m Dh_Z(i) \geq 2r,$$
a contradiction.

In the odd case, we develop the computations only for $m = 2e+1$ and $r = \binom{2e+3}2 +e+1 $, the other  cases 
being covered by Theorem \ref{range}. 

As $r> \binom{2e+3}2$, then $k_{2e+1}(A)= \binom{2e+3}2$, so 
$h_A(2e+1)$ coincides with the dimension of the space of forms in three variables of degree $2e+1$. This implies that the 
evaluation map is injective up to degree $2e+1$, i.e.
$h_A(i)= \binom{i+2}2$ and $Dh_A(i)=i+1$ for $i\leq 2e+1$. It follows that $Dh_Z(i)=i+1=Dh_A(i)$ for $i=0,\dots,2e+1$.
In particular $Dh_Z(2e+1)=2e+2$.

Moreover $Dh_A(2e+2)=r-\binom{2e+3}2= e+1<2e+2.$

It follows from Proposition \ref{prop:nonincr} that  $Dh_Z(i)\leq Dh_Z(i-1)$ for $i\geq 2e+2$.
If $Dh_Z(i)<Dh_Z(i-1)$ for $i=2e+2,\dots,d+1$ until it reaches $0$, then we get $Dh_Z(d+1)= 0$, a contradiction.
Thus there exists $j\geq 2e+2$, $j\leq d$, such that $Dh_Z(j)=Dh_Z(j+1)<Dh_Z(2e+1)$. By Proposition \ref{Dav} we get the contradiction.
\end{proof}

\begin{exa0}\label{sharprn} Even the bounds of Theorem \ref{ranger} are sharp.

The examples are analogous the the ones of Example \ref{sharpn}.

When $d=2m$, $m\geq 5$, we get an example of a form of degree $d$ with one non-redundant decomposition of length  $r=\binom{m+2}2+1$
and one non-redundant decomposition of length $\binom{m+2}2-1$
 by taking a general set of $r$ points and embedding it in a general complete intersection of type $m+1,m+2$.

When $d=4e+1$, $e\geq 2$, we get an example of a form of degree $d$ with one non-redundant decomposition of length  $r=\binom{2e+2}2+e+1$
and one non-redundant decomposition of length $\binom{2e+2}2+e$
 by taking a general set of $r$ points and embedding it in a general complete intersection of type $2e+1,2e+3$.

When $d=4e+3$, $e\geq 2$, we get an example of a form of degree $d$ with one non-redundant decomposition of length  $r=\binom{2e+3}2+e+2$  
and one non-redundant decomposition of length $\binom{2e+3}2 +e$  by taking a general set of $r$ points and embedding it 
in a general complete intersection of type $2e+2,2e+4$.
\end{exa0}

Notice that the even case $d=2m$ of Theorem \ref{ranger} is covered by part (a) of Theorem 1.1 of \cite{Ball19}, while the odd cases extend  the results of \cite{Ball19} and \cite{MourOneto}, for forms in three variables. 
 
A similar situation holds for a general number $n+1$ of variables. We can recover, with the same techniques,
Theorem 1.1 of \cite{Ball19} and Theorem 2.17 of \cite{MourOneto}.

Moreover, by using Theorem  3.6 of  \cite{BigaGerMig94}, one can prove a statement which somehow
extends the previous results.  
Indeed, e.g. in the even case, we show that when $h_A(m-1)$ is not $\ell(A)$, but it is sufficiently closed to $\ell(A)$,
then one can conclude that $T$ is identifiable, thus the rank of $T$ is $\ell(A)$. 

\begin{prop0}\label{MO}
Let $A\subset\Pj^n$ be a non-redundant, non-degenerate set which computes  the form $T$ of degree $d\geq 3$ in $n+1$ variables. 
Put $r=\ell(A)$ and assume $h_A(1)=\min\{n+1, r\}$.
\begin{itemize}
\item If $d=2m$ is even, assume 
$$h_A(m-1)\geq r-\min\left\{\frac{n-1}2,\frac{m-1}2\right\}.$$

\item If $d=2m+1$ is odd, assume $k_m(A)=r$, and 
$$h_A(m-1)\geq r-\min\left\{\frac{n-1}2,\frac{m-1}2\right\}.$$
\end{itemize}
Then $T$ has rank $r$ and it is identifiable.
\end{prop0}
\begin{proof} Let $B$ be another decomposition of $T$, with $\ell(B)\leq r$ and let $Z=A\cup B$. 
By induction on $r$, we can reduce to the case $A\cap B\neq \emptyset$, by applying Remark \ref{redint} just as in the proof
of Theorem \ref{range}. 
Notice indeed that our assumptions on $h_A(m-1)$ are equivalent to say that 
$\sum_{i=m}^\infty Dh_A(i)< \min\{(n-1)/2,(m-1)/2\}$: if the condition holds for $A$, it holds also for any subset $A'$ of $A$.
 Thus assume that $A\cap B=\emptyset$, so that, by Proposition \ref{CBconseq}, $Z$ has the property $\mathit{CB}(d)$.

If $d=2m$, then by assumption $\sum_{i=0}^{m-1}Dh_Z(i)\geq r-\min\{(n-1)/2,(m-1)/2\}$, thus also 
$\sum_{i=m+2}^{d+1}Dh_Z(i)\geq r-\min\{(n-1)/2,(m-1)/2\}$. 

Assume $r<n+1$. Since $m\geq 2$, by Cayley-Bacharach one finds that  
{\small{ $$Dh_Z(m)\le \sum_{i=0}^{d+1} Dh_Z(i) -Dh_Z(0)-Dh_Z(1)-Dh_Z(d)-Dh_Z(d+1) $$
$$ \leq 2r - Dh_A(0)-Dh_A(1) -Dh_Z(0)-Dh_Z(1) \leq 2r- 2(Dh_A(0)+Dh_A(1) )\leq 0.$$ } }
 Thus $Dh_Z(d+1)=0$, a contradiction.

If $r>n+1$, then $h_A(1)=n+1$. We have:
\begin{multline*}
Dh_Z(m)+Dh_Z(m+1)\leq 2r-\sum_{i=0}^{m-1}Dh_Z(i)-\sum_{i=m+2}^{d+1}Dh_Z(i)\\
\leq \min\{n-1, m-1\}<n,m.$$
\end{multline*}
It follows by Proposition \ref{prop:nonincr} that $Dh_Z(i)\geq Dh_Z(i+1)$ for $i\geq m+1$. As in the proof of Theorem 
\ref{range}, if for $i=m+1,\dots,d$ we have $Dh_Z(i)> Dh_Z(i+1)$ until $Dh_Z(i)=0$, then we get $Dh_Z(d+1)=0$,
a contradiction. Thus there exists $j\geq m+1$, $j\leq d$, such that $0<Dh_Z(j)=Dh_Z(j+1)<n$.
By Theorem  3.6 of  \cite{BigaGerMig94} we get that $Z$ is contained in a curve of degree $Dh_Z(j)<n$.
Thus $A$ belongs to a curve of degree $<n$, which cannot span $\Pj^n$, i.e. $h_A(1)=1+Dh_A(1)<n+1$, a contradiction.

The case $d=2m+1$ can be proved similarly.
\end{proof}

Proposition \ref{MO} makes the assumption that $h_A(1)$ is maximal. If this assumption fails, the form $T$ 
 is not {\it coincise}: after  a change of coordinates, $T$ is a form in less than $n+1$ variables. Thus if  $h_A(1)<n+1$, then the number
 $n$ in the bound of the theorem is essentially meaningless for $T$, and the statement would not hold.

\begin{rem0}\label{genRank}
Our methods work also for generic ranks, not only for sub-generic ones. Indeed, in the case of ternary forms, if $ d = 5 $ then Theorem 
\ref{range} provides an alternative proof of Sylvester's Theorem, see also \cite{AngeCMazzon}; if $ d = 4 $ (resp. $ d = 6 $) then, according to 
Theorem \ref{ranger}, a form $ T $ with a sufficiently general decomposition of length $ 6 $ (resp. $ 10 $) has rank $ 6 $ (resp. $ 10 $), 
which is the generic one for this particular class of symmetric tensors. 
\end{rem0}

\section{The identifiability of ternary forms of degree $8$}\label{otto}

As an application of our methods, we can analyze the case of plane optics, i.e. we assume that $ d = 8, n = 2 $, and we fix $ T \in S^{8} \C^{3} $. 

Consider a finite set $ A = \{P_{1},\ldots,P_{r}\} \subset \Pj^{2} $ computing $ T $. 

In this section, we will always assume that $ A $ satisfies the following properties:
\begin{itemize}
\item[(i)] $ A $ is non-redundant,
\item[(ii)] $k_{3}(A) = \min\{10,r\}$,
\item[(iii)] $h_{A}(4) = r $.
\end{itemize}
Notice that these properties can be easily verified by a computer (and our algorithm \ref{algor} will do that).

\indent If $ r \leq 13 $, then, by Theorem \ref{range}, $ T $ is identifiable of rank $ r $. \\
\indent If $ r \in \{14,15\} $, then, by Theorem \ref{ranger}, we can conclude that $ A $ computes the rank of $ T $. In particular, when $ r = 15 $, it has been proved in \cite{RaneSchr00} that the general $ T $ has $ 16 $ decompositions of length $ r $. \smallskip

Therefore we focus on the case $ r = 14 $. \smallskip

In this case we are able to provide a criterion to detect identifiable tensors. In order to do that, we need to introduce the following:

\begin{nota0}\label{dualnot}
From now on, we denote by $ A^{\vee} $ the dual set of $ A $ in $ (\Pj^{2})^{\vee} $, that is $ A^{\vee} = \{P_{1}^{\vee}, \ldots, P_{14}^{\vee}\} $, and by $ J_{A^{\vee}} $ (resp. $ I_{A^{\vee}} $) the ideal sheaf of $ A^{\vee} $ (resp. the ideal defining $ A^{\vee} $).
Moreover, $ (\Pj^{44})^{\vee} $ is the dual space of $ \Pj(S^{8}\C^{3}) \cong \Pj^{44} $ and  $ \mathcal{L} = \langle v_{8}(A) \rangle \cong \Pj^{13} \subset \Pj(S^{8}\C^{3}) $.
\end{nota0}

Since $ A $ satisfies properties (ii) and (iii), then the Hilbert function of $ A $ and its first difference, verify, respectively,
\begin{equation}\label{eq:hDhA}
\begin{tabular}{c|ccccccc}
$j$ & $0$ & $1$ & $2$ &   $3$ &   $4$ &   $5$ &  $\dots$ \\  \hline
$h_{A}(j)$ & $1$ & $3$ &   $6$ &   $10$ &   $14$ &   $14$ &  $\dots$ \cr
$Dh_{A}(j)$ & $1$ & $2$ &   $3$ &   $4$ &   $4$ & $0$ & $\dots$ \cr
\end{tabular}.
\end{equation}

In particular, passing to cohomology in the exact sequence:
$$ 0 \rightarrow J_{A^{\vee}}(s) \rightarrow \mathcal{O}_{\Pj^{2}}(s) \rightarrow \mathcal{O}_{A^{\vee}}(s) \rightarrow 0  $$
for $ s \in \{4, 5, 6\} $, we get that $ I_{A^{\vee}} = (Q,Q_{1}, Q_{2}, Q_{3}, Q_{4}) $, with $ Q \in S^{4} \C^{3} $ and $ Q_{i} \in S^{5} \C^{3} $ for $ i \in \{1,2,3,4\} $. 
Therefore there exist $ q_{i} \in S^{2} \C^{3} $, $ L_{j} \in \C^{3} $ such that the locally free resolution of $ J_{A^{\vee}} $ is 
\begin{equation}\label{eq:idA}
0 \longrightarrow \mathcal{O}_{\Pj^{2}}(-6)^{\oplus 4} \xrightarrow M  \mathcal{O}_{\Pj^{2}}(-4) \oplus \mathcal{O}_{\Pj^{2}}(-5)^{4} \longrightarrow J_{A^{\vee}} \longrightarrow 0.
\end{equation}
where
\begin{equation}\label{M}
M = \begin{pmatrix} q_{1} & q_{2} & q_{3} & q_{4} \\ L_{1} & L_{2} & L_{3} & L_{4} \\ L_{5} & L_{6} & L_{7} & L_{8} \\ L_{9} & L_{10} & L_{11} & L_{12}\\ 
L_{13} & L_{14} & L_{15} & L_{16} \\ \end{pmatrix} 
\end{equation}
is the \emph{Hilbert-Burch matrix} of $ J_{A^{\vee}} $. $ Q,Q_{1}, Q_{2}, Q_{3}, Q_{4} $ coincide, respectively, with $ (-1)^{i} $ times the minor obtained by 
leaving out the $ i $-th row of $M$, $ i \in \{1,2,3,4,5\} $.\\
\indent Now, assume that $ B = \{P''_{1},\ldots,P''_{\ell(B)}\} \subset \Pj^{2} $ is another finite set computing $ T $ such that
\begin{itemize}
\item[(i)] $\ell(B) = 14$;
\item[(ii)] $ B $ is non-redundant
\end{itemize}
and set $ Z = A \cup B \subset \Pj^{2} $. 

\begin{claim0}\label{CB8}
Therefore $ Z $ satisfies $\mathit{CB}(8)$.
\end{claim0}

\begin{proof}
If this is not the case, then, by Proposition \ref{CBconseq}, it holds that $ A \cap B \not= \emptyset $ and so, by means of Remark \ref{redint}, 
there exists $ T_{0} \in S^{8} \C^{3} $ admitting two disjoint non-redundant decompositions $ A' $ and $ B' $, with $ A' \subset A $. If $ \ell(A') < \ell(A) $, 
then $ T_{0} $ has two non-redundant decompositions, with $ \ell(B') \leq \ell(A') \leq 13 $. Since $ A' \subset A $ and $ A $ satisfies properties (ii) and (iii), 
then $ k_{3}(A') = \min \{10, \ell(A')\} $ and $ h_{A'}(4) = \ell(A') \leq 13 $, and so, by Theorem \ref{range}, $ B' $ cannot exist. Thus $ A = A' $. 
Since $ A \cap B' = \emptyset $, from Proposition \ref{CBconseq} we get that $ Z = A \cup B' $ satisfies the property $\mathit{CB}(8)$, which is a contradiction.
\end{proof}

\begin{claim0}\label{cor:uniquecase}
The first difference of the Hilbert function of $ Z $ verifies
\begin{center}\begin{tabular}{c|ccccccccccccc}
$j$ & $0$ & $1$ & $2$ &   $3$ &   $4$ &   $5$ &  $6$ &  $7$ &  $8$ &  $9$ & $ 10 $ & $ \dots $ \\  \hline
$Dh_Z(j)$ & $1$ & $2$ &   $3$ &   $4$ &   $4$ &  $4$ &  $4$ &  $3$ & $2$ & $1$ & $ 0 $ & $ \cdots $ \cr
\end{tabular}\end{center}
Therefore $ A \cap B = \emptyset $, $ \ell(Z) = 28 $ and $ Z^{\vee} $ can be obtained as a complete intersection of type $(4,7)$.
\end{claim0}

\begin{proof}
Notice that, since $ A \subset Z $ and we have \eqref{eq:hDhA}, then $ Dh_{Z}(j) = j+1 $ for $ j \in \{0,1,2,3\} $ and $ Dh_{Z}(4) \geq 4 $. 
Moreover, since $ T $ admits at least two decompositions, then, by Proposition \ref{d+1}, we get that $ Dh_{Z}(9) > 0 $. \\
Now, Claim \ref{CB8} and Theorem \ref{thm:extGKR} imply that:
$$ Dh_{Z}(5)+ \ldots + Dh_{Z}(9) \geq 10 + Dh_{Z}(4) \geq 14 $$
and since:
$$ Dh_{Z}(0)+ \ldots + Dh_{Z}(9) = 10+ Dh_{Z}(4)+Dh_{Z}(5)+ \ldots + Dh_{Z}(9) \leq \ell(Z) \leq 28, $$
then $ Dh_{Z}(5)+ \ldots + Dh_{Z}(9) \leq 14. $ Therefore:
\begin{equation}\label{eq:14}
Dh_{Z}(5)+ \ldots + Dh_{Z}(9) = 14 
\end{equation}
$$ Dh_{Z}(4) = 4. $$
In particular, $ Dh_{Z}(j) = 0 $ for $ j \geq 10 $, $ \ell(Z) = 28 $, $ A \cap B = \emptyset $, and, by Proposition \ref{prop:nonincr}, 
\begin{equation}\label{eq:dis}
4 \geq Dh_{Z}(5) \geq \ldots \geq Dh_{Z}(9). 
\end{equation}
Notice that $ Dh_{Z}(5) \not\in \{1,2\} $. Indeed, if $ Dh_{Z}(5) = 1 $ (resp. $ Dh_{Z}(5) = 2 $), then, by \eqref{eq:dis}, $ Dh_{Z}(5)+ \ldots + 
Dh_{Z}(9) \leq 5 $ (resp. $ Dh_{Z}(5)+ \ldots + Dh_{Z}(9) \leq 10 $), which contradicts \eqref{eq:14}. So, assume that $ Dh_{Z}(5) = 3 $, then, by \eqref{eq:14}, 
$$ Dh_{Z}(6)+ \ldots + Dh_{Z}(9) = 11 $$
so that $ Dh_{Z}(6) = Dh_{Z}(7) = Dh_{Z}(8) = 3 $ and $ Dh_{Z}(9) = 2 $. This fact provides a contradiction thanks to Proposition \ref{Dav}. Thus 
\begin{equation}\label{eq:5}
Dh_{Z}(5) = 4. 
\end{equation}
Notice that $ Dh_{Z}(6) \not\in \{1,2\} $. Thus, suppose that $ Dh_{Z}(6) = 3 $. Then, by \eqref{eq:14} and \eqref{eq:5} it has to be $ Dh_{Z}(7) = 3 $, 
which contradicts Proposition \ref{Dav}, as above. Necessarily, 
\begin{equation}\label{eq:6}
Dh_{Z}(6) = 4.
\end{equation}
Therefore, by \eqref{eq:14}, \eqref{eq:5} and \eqref{eq:6}, 
$$ Dh_{Z}(7)+ Dh_{Z}(8) + Dh_{Z}(9) = 6. $$
If $ Dh_{Z}(7)= Dh_{Z}(8)= Dh_{Z}(9) = 2 $, then we get again a contradiction by Proposition \ref{Dav}. Thus it has to be
$$ Dh_{Z}(7) = 3, Dh_{Z}(8) = 2 , Dh_{Z}(9) = 1, $$
as desired.
In particular, by Theorem \ref{thm:Davis}, $ Z^{\vee} $ is contained in a plane quartic. Moreover, passing to cohomology in the exact sequence
$$ 0 \rightarrow J_{Z^{\vee}}(s) \rightarrow \mathcal{O}_{\Pj^{2}}(s) \rightarrow \mathcal{O}_{Z^{\vee}}(s) \rightarrow 0  $$
for $ s \in \{4, 7, 11\} $, we get that $ Z^{\vee} $ is contained in a unique quartic $Q$, and there exists a septic containing $Z$ and not containing $Q$. 
Since, $ Z $ satisfies $ CB(8) $ and the Hilbert function of $Z$ is the same as the Hilbert function of a complete intersection of type $ (4,7) $, then, 
by the Main Theorem of \cite{Davis84}, $ Z^{\vee} $ \emph{is} a complete intersection of type $ (4,7) $, which allows us to conclude the proof.
\end{proof}

As a consequence of Claim \ref{cor:uniquecase}, $ I_{Z^{\vee}} = (Q,S) $, where $ Q \in S^{4}\C^{3} $ and $ S \in S^{7}\C^{3} $. In particular, 
$ Q \in H^{0}(J_{A^{\vee}}(4)) $ and $ S \in H^{0}(J_{A^{\vee}}(7)) $. By applying Proposition 5.2.10 of \cite{Migliore} (Mapping cone) to the commutative diagram 
$$\begin{CD}
0 @>>> \mathcal{O}_{\Pj^{2}}(-6)^{\oplus 4} @> M  >> \mathcal{O}_{\Pj^{2}}(-4)\oplus \mathcal{O}_{\Pj^{2}}(-5)^{\oplus 4} @>>> J_{A^{\vee}} @>>> 0\\
@. @AA M_{1} A @AA M_{2}   A @AAA @.  \\
0 @> >>\mathcal{O}_{\Pj^{2}}(-11) @>>{\it \begin{pmatrix} -S \\ Q \end{pmatrix}} >\mathcal{O}_{\Pj^{2}}(-4)\oplus \mathcal{O}_{\Pj^{2}}(-7) @>>> J_{Z^{\vee}} @>>> 0
\end{CD}$$
where $ M $ satisfies \eqref{M} and
$$ M_{1} = \begin{pmatrix} Q'_{1} \\ Q'_{2} \\ Q'_{3} \\ Q'_{4} \end{pmatrix}, \, M_{2} ={\it \begin{pmatrix} a &  0 \\  0 & q'_{1} \\ 0 & q'_{2} \\ 0 & q'_{3} \\ 0 & q'_{4} \\ 
\end{pmatrix}}$$
with $ Q'_{i} \in S^{5} \C^{3}, a \in \C, q'_{j} \in S^{2}\C^{3} $, for $ i,j \in \{1,2,3,4\} $, we get that $ J_{B^{\vee}} $ has a locally free resolution of the form 
\begin{equation}\label{eq:idB}
0 \rightarrow \mathcal{O}_{\Pj^{2}}(-6)^{\oplus 4}  \xrightarrow {SM} \mathcal{O}_{\Pj^{2}}(-4) \oplus \mathcal{O}_{\Pj^{2}}(-5)^{\oplus 4} \rightarrow J_{B^{\vee}} 
\rightarrow 0
\end{equation}
where
$$ SM = \begin{pmatrix}q'_{1} & q'_{2} & q'_{3} & q'_{4} \\  L_{1} & L_{5} & L_{9}  & L_{13} \\ L_{2} & L_{6} & L_{10} & L_{14} \\ L_{3} & L_{7} & L_{11} & L_{15} \\ 
L_{4} & L_{8} & L_{12} & L_{16}  \\ \end{pmatrix}. $$ \\
Notice that the lower part of the matrix $SM$ is the transpose of the lower part of the matrix $M$.

Since $ \dim (I_{A^{\vee}})_4 = 1 $, finite sets $ B^{\vee} $ obtained from $ A^{\vee} $ in a complete intersection of type $ (4,7) $ 
are parameterized by the space of septics in $I_A$, modulo those septics which are multiples of  $Q$, i.e. by the space:
$$\frac{ (I_{A^\vee})_7}{S^2\C^3\cdot (I_{A^\vee})_4},$$ 
which is a linear space of projective dimension 
$$ 11 = \dim(S^{7}\C^{3})-\ell(A)-\dim(S^{3}\C^{3})-1 = 36-14-10-1. $$
Therefore we can identify such a set $ B^{\vee} $ with an element 
 in $  H^{0}(J_{A^{\vee}}(7)) /(S^3(\C^3)\otimes H^{0}(J_{A^{\vee}}(4))).$ Thus any $S\in H^{0}(J_{A^{\vee}}(7))$ which is  not a multiple of
the quartic $Q$ determines a set $B^{\vee}$, and we will denote it by $ B(S)^{\vee} $. In order to get all such finite sets $ B(S)^{\vee} $, 
it suffices to focus on the matrix $ SM $ of \eqref{eq:idB}, where the $ L_{j}'s $ are fixed while the $ q'_{j} $'s depend on $ 24 $ parameters, let us say 
$$ q_{j}' = a_{0+6j}x_{0}^2+2a_{1+6j}x_{0}x_{1}+2a_{2+6j}x_{0}x_{2}+a_{3+6j}^2x_{1}^2+2a_{4+6j}x_{1}x_{2}+a_{5+6j}x_{2}^2 $$
with $ j \in \{1,2,3,4\}. $ By applying elementary rows operations to $ SM $, we can assume, without loss of generality, that $ q'_{1} = q'_{3} = 0 $, 
so that the parameters reduce to $ 12 $. More in detail, consider the polynomial system
\begin{equation}\label{eq:system}
\left\{\begin{array}{rcl} \ell_1 L_1+\ell_2 L_2+\ell_3 L_3+\ell_4 L_4& =& q'_1\\
\ell_1 L_9+\ell_2 L_{10}+\ell_3 L_{11}+\ell_4 L_{12}& =& q'_3\\
\end{array}\right.
\end{equation} 
where $ \ell_1,\ell_2,\ell_3,\ell_4 \in \C^3 $ are unknown. Assume that $ \ell_{i} = l_{i0}x_{0} + l_{i1}x_{1} + l_{i2}x_{2}$, for $ i \in \{0,1,2\} $, 
with $ \ell_{ij} \in \C $, and recall that $ L_{j} = (\partial_{0} L_{j})x_{0} + (\partial_{1} L_{j})x_{1} +(\partial_{2} L_{j})x_{2}$, for any$ j $, where 
$ \partial_{h}L_{j} \in \C $ denotes the partial derivative of $ L_{j} $ with respect to $ x_{h} $. From \eqref{eq:system}, we get a linear system with 
$ 12 $ equations in $ 12 $ unknowns (the coefficients of the $ \ell_{i}'s $). Let $ C $ be the $ 12 \times 12 $ matrix associated to the system.
Direct computations show that:
$$ C = \begin{pmatrix}\partial_{0}L_{1} & 0 & 0 & \cdots &  \cdots  &  \partial_{0}L_{4} & 0 & 0 \\ 
\partial_{1}L_{1} & \partial_{0}L_{1} & 0 & \cdots & \cdots &  \partial_{1}L_{4} & \partial_{0}L_{4} & 0\\
\partial_{2}L_{1} & 0 & \partial_{0}L_{1} & \cdots & \cdots &  \partial_{2}L_{4} & 0 & \partial_{0}L_{4}\\
0 & \partial_{1}L_{1} & 0 & \cdots &  \cdots &  0 & \partial_{1}L_{4} & 0\\
0 & \partial_{2}L_{1} & \partial_{1}L_{1} & \cdots & \cdots &  0 & \partial_{2}L_{4} & \partial_{1}L_{4}\\
0 & 0 & \partial_{2}L_{1} & \cdots & \cdots & 0 & 0 & \partial_{2}L_{4}\\
\partial_{0}L_{9} & 0 & 0 & \cdots &  \cdots  &  \partial_{0}L_{12} & 0 & 0\\ 
\partial_{1}L_{9} & \partial_{0}L_{9} & 0 & \cdots &  \cdots  &  \partial_{1}L_{12} & \partial_{0}L_{12} & 0\\
\partial_{2}L_{9} & 0 & \partial_{0}L_{9} & \cdots &  \cdots &  \partial_{2}L_{12} & 0 & \partial_{0}L_{12}\\
0 & \partial_{1}L_{9} & 0 & \cdots &  \cdots &  0 & \partial_{1}L_{12} & 0\\
0 & \partial_{2}L_{9} & \partial_{1}L_{9} & \cdots &  \cdots  &  0 & \partial_{2}L_{12} & \partial_{1}L_{12}\\
0 & 0 & \partial_{2}L_{9} & \cdots & \cdots & 0 & 0 & \partial_{2}L_{12}\\ 
\end{pmatrix}.$$
One computes, either by hand or with \cite{Macaulay2} that $ C $ has full rank $ 12 $, for a generic choice of $ L_{1}, L_{2}, L_{3}, L_{4} $ (one random example suffices). 
By Kramer's theorem, \eqref{eq:system} admits a unique solution.

Let $ A,B(S) \subset \Pj^{2}$ be as above and let 
$$ \Pj_{A} = \Pj(H^{0}(J_{A^{\vee}}(8))) \cong \Pj^{30},$$ 
so that $\Pj_{A}$ determines a linear space of dimension $44-14 = 30$ inside the dual space $\Pj(S^8 \C^3)^{\vee} \cong (\Pj^{44})^{\vee}$, 
which has been introduced in Notation \ref{dualnot}.

Define similarly
$$ \Pj_{B(S)} = \Pj(H^{0}(J_{B(S)^{\vee}}(8))) \cong \Pj^{30} \subset \Pj(S^8 \C^3)^{\vee} \cong (\Pj^{44})^{\vee}. $$
 By construction, 
$$ \dim (\Pj(H^{0}(J_{A^{\vee}\cup B(S)^{\vee}}(8)))) = \dim (S^{4} \C^{3})+\dim (\C^{3}) - 1 = 17. $$
Therefore, by Grassmann's formula for projective spaces, 
$$ \dim (\Pj(H^{0}(J_{A^{\vee}}(8)) +  H^{0}(J_{B(S)^{\vee}}(8)))) = 30+30-17 = 43  $$
that is $ \Pj(H^{0}(J_{A^{\vee}}(8))+ H^{0}(J_{B(S)^{\vee}}(8))) \subset (\Pj^{44})^{\vee} $ is a hyperplane. By duality, it corresponds to a point 
$\Pj(H^{0}(J_{A^{\vee}}(8)) +  H^{0}(J_{B(S)^{\vee}}(8)))^{\vee}$ of the subspace $\mathcal{L} \subset \Pj^{44} $ introduced in Notation \ref{dualnot}, 
admitting \emph{at least two} decompositions of length $ 14 $, $A$ and $B(S)$, thus it corresponds to the fixed plane optic $ T $. \\
\indent We define in this way a (rational) map
$$ f : \Pj^{11} \dasharrow \mathcal{L} $$
\begin{equation}\label{eq:f}
f(S) = \Pj(H^{0}(J_{A^{\vee}}(8))+ H^{0}(J_{B(S)^{\vee}}(8)))^{\vee}. 
\end{equation}

\begin{claim0}\label{thm:bir}
When $A$ satisfies the assumptions at the beginning of this section, the map $ f : \Pj^{11} \dasharrow \mathcal{L} $ defined in \eqref{eq:f} is birational.
\end{claim0}

\begin{proof}
It suffices to show that for some $ P \in \imm(f) \subset \mathcal{L} $ the set $ f^{-1}(P) $ is finite and has degree $ 1 $. \\
We prove this fact via a computational approach in Macaulay2 \cite{Macaulay2} (over a finite field, but then the proof holds
also over $\C$), see the ancillary file available at \texttt{https://arXiv.org/src/1901.01796v4/anc/optics.txt}. In particular, we fix a finite set 
$ A = \{P_1,\ldots,P_{14}\} \subset \Pj^2 $ whose elements have random coefficients. We construct the Hilbert-Burch matrix of $ J_{A^{\vee}} $ 
and we fix an element $ S \in H^0 (J_{A^{\vee}}(7))$, not multiple of the quartic $Q$.
This is equivalent to a choice of $4$ conics $ q'_{1}, q'_{2}, q'_{3}, q'_{4} $ (with $ q'_{1} = q'_{3} = 0 $, $ q'_{2},q'_{4} \neq 0 $) and so of a residual 
set $ B(S)^{\vee} $ whose ideal sheaf admits a free resolution as in \eqref{eq:idB}. By means of \eqref{eq:f}, we compute $ f(S) $ and we pose 
$ P = f(S) $. Let $ (p_{0}, \ldots, p_{44}) $ be a representative vector for the point $ P $. \\
In order to get $ f^{-1}(P) $, in the first row of the Hilbert-Burch matrix $ SM $ of $ J_{B(S)^{\vee}} $ we change $ q'_{j} $ with $ q_{j}' = a_{0+6j}x_{0}^2+
2a_{1+6j}x_{0}x_{1}+2a_{2+6j}x_{0}x_{2}+a_{3+6j}^2x_{1}^2+2a_{4+6j}x_{1}x_{2}+a_{5+6j}x_{2}^2   $, for $ j \in \{2,4\} $ and we consider the 
$ 45\times 44 $ matrix $ M$Fix$'' $ whose columns are a set of generators  for $ H^0(J_{A^{\vee}}(8))+ H^0(J_{B(S)^{\vee}}(8)) $. 
 Notice that $ M$Fix$''$ is divided in $ 2 $ blocks: the first $ 31 $ columns have integer entries while in the last $ 13 $ the entries depend linearly on the 
 $ 12 $ parameters $ a_{6}, \ldots, a_{11}, a_{18}, \ldots, a_{23}  $. Let us say  $ M$Fix$'' = A_{1}|A_{2} $. Therefore 
$$ f^{-1}(P) = \{(a_{6}, \ldots, a_{11}, a_{18},\ldots, a_{23}) \in A^{12} \, | \, (p_{0}, \ldots, p_{44})\cdot MFix'' = 0_{1 \times 45}\}, $$
where $ A^{12} $ denotes the affine space of dimension $ 12 $.
Since $ (p_{0}, \ldots, p_{44})\cdot A_{1} = 0_{1 \times 31} $ provide trivial conditions, then 
\begin{equation}\label{eq:linsys}
f^{-1}(P) = \{(a_{6}, \ldots, a_{11}, a_{18},\ldots, a_{23}) \in A^{12} \, | \, (p_{0}, \ldots, p_{44})\cdot A_{2} = 0_{1 \times 13}\}. 
\end{equation}
The $ 13 \times 12 $ matrix associated to the linear system appearing in \eqref{eq:linsys} has rank $ 11 $. Then, by Kramer's  theorem, the 
affine dimension of $ f^{-1}(P) $ is $ 1 $, which allows us to conclude the proof.
\end{proof}

Claim \ref{thm:bir} implies the following:

\begin{claim0}
If $ T \in S^{8}\C^{3} $ of rank $ 14 $ is a general point in the image of $f$, i.e. a general unidentifiable optic 
of rank $14$, then there are \emph{exactly two} finite sets computing the rank of $ T $. 
\end{claim0}

As a consequence we get the following:

\begin{claim0}
$ \mathcal{L} $ contains a variety of projective dimension $ 11 $, whose general points consist of forms in $ S^{8}\C^{3} $ of rank $ 14 $ that 
admit two finite sets computing the rank. 
\end{claim0}

Now we are able to explain a relevant consequence of our analysis:

\begin{rem0}\label{rem:crit} From the construction outlined above, one can develop a criterion that, given $ T \in S^{8}\C^{3} $ of rank $ 14 $ admitting a 
non-redundant finite set $ A = \{P_{1}, \ldots, P_{14}\} \subset \Pj^{2} $ computing it with $ k_{3}(A) = 10 $ and $ h_{A}(4) = 14 $, 
establishes the uniqueness of such an $ A $, i.e. the identifiability of $T$. 

Indeed, if the rank of $ 13 \times 12 $ matrix of the linear system in \eqref{eq:linsys} is $ 12 $, then $ A $ is unique. 
\end{rem0}

In what follows we describe the algorithm based on the criterion introduced in Remark \ref{rem:crit}.

\subsection{The algorithm} \label{algor} Consider a finite set $ A = \{P_{1}, \ldots, P_{14}\} \subset \Pj^{2} $ in the form of a collection of points 
$ A^{\vee} = \{P_i^{\vee} = [\vect{v}_i] \}_{i=1}^{14} \subset (\Pj^2)^{\vee}$ and a ternary form $ T $ of degree $ 8 $ in the linear span of $ v_{8}(A) $. 
This means that for any choice of representatives (coordinates) $T_1,\dots, T_{14}$ for the projective points $v_8(P_1),\dots,v_8(P_{14})$ respectively, 
we have 
\[
 T = \sum_{i=1}^{14} \lambda_{i}T_i 
\]
for certain $ \lambda_{1}, \ldots \lambda_{14}  \in \C $, moreover $T$ represents the projective point $ [(p_{0}, \ldots, p_{44})]^\vee$.
According to Theorem \ref{ranger} we can perform the next tests for verifying that $T$ has rank $ 14 $:
 \begin{enumerate}
  \item[1)] \emph{non-redundanty test}: check that $\dim \langle v_8(\vect{v}_1), \ldots, v_8(\vect{v}_{14}) \rangle = 14$;
  \item[2)] \emph{fourth Hilbert function test}: check that $h_{4}(A) = 14 $. 
 \end{enumerate}
 If all these tests are successful, then $T$ is of rank $14$. \\
With the notation introduced in the proof of Claim \ref{thm:bir}, if, in addition, the following tests provide positive answers:
 \begin{enumerate}
\item[3)] \emph{third Kruskal's rank test}: check that $k_{3}(A) = 10 $,
\item[4)] check that the $ 13 \times 12 $ matrix of the linear system $(p_{0}, \ldots, p_{44})\cdot A_{2} = 0_{1 \times 13}$ has rank $ 12 $,
 \end{enumerate}
then $ f^{-1}(T) $ is empty and so $ T $ is identifiable. 

The algorithm has been implemented in Macaulay2, over the finite field $ \Z_{31991} $. For more details, see the ancillary file 
available at \texttt{https://arXiv.org/src/1901.0\\1796v4/anc/optics.txt}.

This new criterion is effective in the sense of \cite{COttVan17b}. Indeed, ternary forms computed by $ 14 $ summands are generically identifiable 
\cite{COttVan17a}, and it is easy to verify that the conditions in tests $1)$, $2)$, $3)$ and $ 4) $ are not satisfied precisely on a Zariski-closed 
strict sub-variety of the $14$-secant variety of $v_8(\Pj^2)$. \\

In the next subsection, we present some examples of identifiable and unidentifiable ternary forms of degree $ 8 $ and rank $ 14 $.

\subsection{Examples}\label{exok}
In Macaulay2, we generated a random collection of $ 14 $ points $ A^{\vee} = \{P_{i}^{\vee}= [\vect{v}_i] \}_{i=1}^{14} $, where 
\[
\begin{bmatrix}
 \vect{v}_i
\end{bmatrix}_{i=1}^{14} = {\small{
\begin{bmatrix}
42 & -4 & 17 \\
-50 & -36 &-28  \\
39 & -16 & 37 \\
9 & -6 & -22 \\
-15 & -32 & -19\\
-22 & 31 & 45 \\
50 & -32 & -8 \\
45 & -38 & -31 \\
-29 & 31 & -9 \\
-39 & 24 & 32 \\
30 & -42 & -4 \\
19 & -50 & 4 \\
-38 & -41 & -2 \\
2 & 15 & 24
\end{bmatrix}.}}
\]
The non-redundanty test shows that $\dim \langle v_8(A) \rangle = \operatorname{rank}([v_{8}(\vect{v}_i)]_{i=1}^{14}) = 14$, as desired. We then compute $ h_{A}(4)= \operatorname{rank}([v_{4}(\vect{v}_i)]_{i=1}^{14}) $, getting $ 14 $ as required. Notice that these two conditions are satisfied for any $ T \in \langle v_{8}(A) \rangle $. Therefore, any $ T $ computed by $ A $ has rank $ 14 $. \\
Finally we compute the rank of all $ 1001 $ subsets of $ 10 $ columns of $ [v_{3}(\vect{v}_i)]_{i=1}^{14} $. They are all of rank $ 10 $ and so $ k_{3}(A) = 10 $. As for the previous tests, this condition holds for any $ T \in \langle v_{8}(A) \rangle $.\\ Therefore, the identifiability of $ T \in \langle v_{8}(A) \rangle $ depends on the choice of the coefficients $\lambda_i$'s that express $T$ as a linear combination of the points $v_8(P_i)$'s. \\

\paragraph{\textit{An identifiable case}}  Let 
$$ T_{1} = \sum_{i=1}^{14} v_{8}(P_{i}) = $$
{\small{$$ = [-4160x_{0}^8+10086x_{0}^7x_{1}-10592x_{0}^6x_{1}^2-13805x_{0}^5x_{1}^3-5415x_{0}^4x_{1}^4-728x_{0}^3x_{1}^5-10682x_{0}^2x_{1}^6+$$
$$+11924x_{0}x_{1}^7+11680x_{1}^8-10568x_{0}^7x_{2}+2172x_{0}^6x_{1}x_{2}+4949x_{0}^5x_{1}^2x_{2}-12129x_{0}^4x_{1}^3x_{2}+$$ $$+10744x_{0}^3x_{1}^4x_{2}+2672x_{0}^2x_{1}^5x_{2}-12873x_{0}x_{1}^6x_{2}-1107x_{1}^7x_{2}-9188x_{0}^6x_{2}^2+9276x_{0}^5x_{1}x_{2}^2+$$ $$+732x_{0}^4x_{1}^2x_{2}^2+11721x_{0}^3x_{1}^3x_{2}^2-13726x_{0}^2x_{1}^4x_{2}^2+3431x_{0}x_{1}^5x_{2}^2-8124x_{1}^6x_{2}^2+12437x_{0}^5x_{2}^3+$$
$$+15504x_{0}^4x_{1}x_{2}^3+9356x_{0}^3x_{1}^2x_{2}^3-14840x_{0}^2x_{1}^3x_{2}^3-4473x_{0}x_{1}^4x_{2}^3+2175x_{1}^5x_{2}^3-12329x_{0}^4x_{2}^4+$$
$$-1390x_{0}^3x_{1}x_{2}^4+6775x_{0}^2x_{1}^2x_{2}^4-2372x_{0}x_{1}^3x_{2}^4-9493x_{1}^4x_{2}^4-7958x_{0}^3x_{2}^5-13661x_{0}^2x_{1}x_{2}^5+$$
$$-11117x_{0}x_{1}^2x_{2}^5+3342x_{1}^3x_{2}^5-5685x_{0}^2x_{2}^6-9054x_{0}x_{1}x_{2}^6+1829x_{1}^2x_{2}^6-1350x_{0}x_{2}^7+7453x_{1}x_{2}^7+12146x_{2}^8] .$$}}

\noindent Since tests 1) and 2) are successful, then $ T_{1} $ has rank $ 14 $. Moreover, test 3) provides positive answer and in this case the 
$ 13 \times 12 $ matrix of the linear system $(p_{0}, \ldots, p_{44})\cdot A_{2} = 0_{1 \times 13}$ has rank $ 12 $. Therefore we can conclude 
that $ A $ is the unique non-redundant finite set of length $ 14 $ computing $ T_{1} $.\\

\paragraph{\textit{An unidentifiable case}}  Let {\small{$$ (\lambda_{1}, \ldots, \lambda_{14}) = (-6395,-1019,2227,13599,-2136,-1329,5500,$$
$$ \quad\quad\quad\quad\quad\quad\quad -4082,7252,-2038,13457,8366,8750,-10807) $$}}
and let
$$ T_{2} = \sum_{i=1}^{14}\lambda_{i} v_{8}(P_{i}) = $$
{\small{$$ = [14990x_{0}^8+748x_{0}^7x_{1}+1813x_{0}^6x_{1}^2-1788x_{0}^5x_{1}^3-8326x_{0}^4x_{1}^4+3614x_{0}^3x_{1}^5-6672x_{0}^2x_{1}^6-6515x_{0}x_{1}^7+$$
$$-5729x_{1}^8+8254x_{0}^7x_{2}-1824x_{0}^6x_{1}x_{2}+1630x_{0}^5x_{1}^2x_{2}-5694x_{0}^4x_{1}^3x_{2}-2192x_{0}^3x_{1}^4x_{2}+12142x_{0}^2x_{1}^5x_{2}+$$
$$-10283x_{0}x_{1}^6x_{2}-6291x_{1}^7x_{2}+13369x_{0}^6x_{2}^2-5192x_{0}^5x_{1}x_{2}^2+11695x_{0}^4x_{1}^2x_{2}^2+8920x_{0}^3x_{1}^3x_{2}^2+$$
$$+11932x_{0}^2x_{1}^4x_{2}^2+10224x_{0}x_{1}^5x_{2}^2+15877x_{1}^6x_{2}^2+6491x_{0}^5x_{2}^3-1780x_{0}^4x_{1}x_{2}^3+9943x_{0}^3x_{1}^2x_{2}^3+$$
$$-109x_{0}^2x_{1}^3x_{2}^3-9947x_{0}x_{1}4x_{2}^3+8699x_{1}^5x_{2}^3-12334x_{0}^4x_{2}^4-14722x_{0}^3x_{1}x_{2}^4+5584x_{0}^2x_{1}^2x_{2}^4+$$
$$+14422x_{0}x_{1}^3x_{2}^4-11037x_{1}^4x_{2}^4+15296x_{0}^3x_{2}^5+15632x_{0}^2x_{1}x_{2}^5-909x_{0}x_{1}^2x_{2}^5-11303x_{1}^3x_{2}^5-12198x_{0}^2x_{2}^6+$$
$$+3575x_{0}x_{1}x_{2}^6+4010x_{1}^2x_{2}^6-12257x_{0}x_{2}^7+11144x_{1}x_{2}^7+x_{2}^8]. $$}}

\noindent As in the previous case, tests 1), 2) and 3) are successful for $ T_{2} $. Notice that $ T_{2}^{\vee} \in \imm(f) $. Indeed, $ T_{2} = f(S) $, 
where $ S $ is the plane septic defined by the vanishing of the determinant of the $ 5\times 5 $ matrix obtained by adding the row 
$ (0,x_{0}^2-x_{0}x_{1}-x_{0}x_{2}+x_{1}^2-x_{1}x_{2}+x_{2}^2, 0, x_{0}^2+x_{0}x_{1}-x_{0}x_{2}+x_{1}^2+x_{1}x_{2}+x_{2}^2) $ to the transpose 
of the Hilbert-Burch matrix $ M $ of $ J_{A^{\vee}} $. In particular $ T_{2} $ is computed by two non-redundant finite sets of length $ 14 $, 
$ A $ and $ B(S) $: the latter is the residual set of the former in the complete intersection of type $ (4,7) $ given by the unique plane 
quartic $ Q $ passing through $ A^{\vee} $ and the plane septic $ S $. According to our theory, test 4) must fail, since $ T_{2} $ is unidentifiable. 
Performing the computation, we find that the $ 13 \times 12 $ matrix of the linear system $(p_{0}, \ldots, p_{44})\cdot A_{2} = 0_{1 \times 13}$ 
has rank $ 11$, one less than expected. In particular, it follows that $ f^{-1}(T_{2}) $ consists of a singleton and $ T_{2} $ 
is computed by exactly two finite sets of length $ 14 $, $ A $ and $ B(S) $.

\begin{rem0}
The second decomposition $ B(S) $ of $ T_{2} $ can be recovered by means of standard numerical methods. Indeed, generically the points of 
$ A^{\vee} \cup B(S)^{\vee} $ are contained in an affine chart of $ \Pj^{2} $ and so, according to the theory of resultants developed in \cite{CLOS}, 
the coefficients of the $ 28 $ linear factors of $ Res_{1,4,7}(F_{0},Q,S) $, where $ F_{0} = u_{0}x_{0}+u_{1}x_{1}+u_{2}x_{2} $ and $  u_{0}, u_{1}, u_{2} $ 
are independent variables, provide the points of $ A^{\vee} \cup B(S)^{\vee} $. Since $ A^{\vee} $ is known, this method yields $ B(S)^{\vee} $. 
\end{rem0}

\smallskip

\begin{rem0}
Similar phenomena occur for higher degrees and for ranks that approximate the generic one. For example, consider the case of plane curves of 
degree $9$ admitting a non-redundant decomposition of length $ 18 $. In this setting, depending of the coefficient of the expression of the form 
$T$ in terms of the decomposition, it may happen that $T$ is identifiable or unidentifiable. 

Moreover, a new phenomenon occurs. Even if the decomposition $A$ is general (\emph{and non-redundant}), 
the rank of the form under investigation may be lower than $ 18 $. In other words, there might exist another decomposition
$B$ of $T$ with only $17$ points.

These examples will be the object of a forthcoming paper.
\end{rem0}

\begin{rem0}\label{nonsing}
With the notation of Example \ref{exok},  forms $T$ of degree $8$ corresponds to points of the secant variety $\sigma_{14}$ to the Veronese variety $v_8(\Pj^2)$.
A non trivial geometric question concerns the description of the singularities of secant varieties. We refer to the paper 
of Han \cite{Han18} for an account on the problem.

Our construction determines the existence of points, in the span of a general set of $14$ points $v_8(A)$ in $v_8(\Pj^2)$,
which are singular, and even non-normal, for the secant variety $\sigma_{14}(v_8(\Pj^2))$.

Here is the reason: the abstract secant variety  $\Sigma_{14}$ (see section 4 of \cite{C19}) is smooth at the point $(T_2,v_8(A))\in \Pj^{44}\times (v_8(\Pj^2))^{14}$,
where $T$ is  the form $T_2$ defined in Example \ref{exok}, 
because it is locally a $\Pj^{13}$-bundle over $v_8(\Pj^2)$, around $(T,v_8(A))$ (for $v_8(A)$ is linearly independent). The form $T$
has exactly two different decompositions $A,B$, thus there are two points of $\Sigma_{14}$,  $(T,v_8(A))$ and $(T,v_8(B))$, which map to $T$.

By the Zariski Main Theorem  (see \cite{Hartshorne}, Corollary 11.4), then $T$ is non-normal in $\sigma_{14}(v_8(\Pj^2)$, unless there exists an infinite family
of points $\eta_t\in \Sigma_{14}$, containing $(T,v_8(A))$) whose points map to $T$. Since $T$ has only two decompositions, the general element $\eta$ 
of the family must be in the closure of $\Sigma_{14}$, and not of type $(T_0,v_8(A_0))$, where $A_0$ is a decomposition of $T_0$.
Since all the points of $\Sigma_{14}$ around $(T,v_8(A))$ are of type  $(T_0,v_8(A_0))$ with $A_0$ decomposition of $T_0$, the family  $\eta_t$ cannot exists.
Thus $T$ is a non-normal point of $\sigma_{14}(v_8(\Pj^2))$.

Similar singular points can be constructed for higher values of the degree $d$, from the examples described in Example \ref{sharpn}.
\end{rem0}

\begin{rem0}
Let us finish with a discussion of how our analysis can be used towards the solution of the following theoretical problem.

{\it Determine an algorithm that can guarantee that a form $T$ represents a smooth point of the secant variety $\sigma_r(X)$,
when $X$ is a Veronese image of $\Pj^2$.}

Consider the case where $X=v_8(\Pj^2)$, $r=14$, and assume one knows a form $T$ which
has a unique decomposition $A$ of length $14$. 

Is then $T$ smooth in $\sigma_{14}(v_8(\Pj^2))$?

Since $A$ is necessarily linearly independent, then  $(T,v_8(A))$ is a smooth point of the corresponding abstract
secant variety $\Sigma_{14}$. Namely $\Sigma_{14}$ is locally a projective bundle over the locus in $X^{r}$
determined by $r$-tuples of linearly independent points.

Thus, $T$ is non-singular in $\sigma_{14}(v_8(\Pj^2))$, unless one of the following two situations holds:
\begin{itemize}
\item In the projection of $\Sigma_{14}$ to $\Pj^{44}$, the tangent space to $\Sigma_{14}$ at $(T,v_8(A))$ drops rank.
As discussed in Section 6 of \cite{AngeCVan18}, this can be excluded if the \emph{Terracini's test holds}: the
dimension of the span $\Theta$ of the tangent spaces to $v_8(\Pj^2)$ at the points of $v_8(A)$ has (the expected) value $41$.
Since $\Theta$ corresponds, in the space of ternary forms of degree $8$, to the degree $8$ part of the ideal spanned
by $L_1^7,\dots, L_{14}^7$, where $L_1,\dots,L_{14}$ are the points of $A$, then a simple computation
on the coefficient matrix of the products $L_i\elena{^7}x_j$, $j=0,1,2$, can exclude this possibility.
\item There exists another point $\eta\in \Sigma_{14}$ which is mapped to $T$ by the projection $\Sigma_{14}\to \Pj^{44}$.  
Here the situation is more involved. Our algorithm above can exclude the existence of such a point $\eta$ of the form $(T,v_8(B))$, 
where $B$ is another decomposition of $T$. 
\end{itemize}

Thus, if we restrict the natural map $\Sigma_{14}\to\sigma_{14}(v_8(\Pj^2))$ to the
strict abstract secant locus $\Sigma'$, defined as the (constructible) set of pairs $(T,v_8(A))$, where $A$ is reduced, $v_8(A)$ is linearly independent
and $T$ belongs to the span of $v_8(A)$, then the algorithm \ref{algor} can guarantee that $T$ represents a non-singular point of the image,
in the sense of complex differential geometry, because it can exclude the existence of $\eta\in\Sigma'$.

On the other hand, the previous method will not be sufficient to guarantee the smoothness of $T$ as a point of $\sigma_{14}(v_8(\Pj^2))$.
Indeed the abstract secant varieties can have (limit) points $\eta\notin\Sigma'$ corresponding to pairs
$(T,v_8(B))$ where $B$ is \emph{not} a decomposition of $T$. 
This can happen e.g. when $B$ is a non-reduced subscheme of length
$14$ of $\Pj^2$ (a \emph{cactus} decomposition). Or, even worse, it may happen that $B$ is a limit, reduced but linearly
dependent set, in which case $T$ is not forced to stay in the span of $v_8(B)$.  

In all the latter cases, a more detailed theoretical analysis to exclude the existence of $\eta$ seems necessary in order
to guarantee the smoothness of $T$.

\end{rem0}

\bibliographystyle{amsplain}
\bibliography{biblioLuca}

\end{document}